# SPIKE AND SLAB VARIABLE SELECTION: FREQUENTIST AND BAYESIAN STRATEGIES

By Hemant Ishwaran[1] and J. Sunil Rao[2]

*Cleveland Clinic Foundation and Case Western Reserve University*

Variable selection in the linear regression model takes many apparent faces from both frequentist and Bayesian standpoints. In this paper we introduce a variable selection method referred to as a rescaled spike and slab model. We study the importance of prior hierarchical specifications and draw connections to frequentist generalized ridge regression estimation. Specifically, we study the usefulness of continuous bimodal priors to model hypervariance parameters, and the effect scaling has on the posterior mean through its relationship to penalization. Several model selection strategies, some frequentist and some Bayesian in nature, are developed and studied theoretically. We demonstrate the importance of selective shrinkage for effective variable selection in terms of risk misclassification, and show this is achieved using the posterior from a rescaled spike and slab model. We also show how to verify a procedure's ability to reduce model uncertainty in finite samples using a specialized forward selection strategy. Using this tool, we illustrate the effectiveness of rescaled spike and slab models in reducing model uncertainty.

**1. Introduction.** We consider the long-standing problem of selecting variables in a linear regression model. That is, given $n$ independent responses $Y_i$, with corresponding $K$-dimensional covariates $\mathbf{x}_i = (x_{i,1}, \ldots, x_{i,K})^t$, the problem is to find the subset of nonzero covariate parameters from $\boldsymbol{\beta} = (\beta_1, \ldots, \beta_K)^t$, where the model is assumed to be

(1) $\quad Y_i = \alpha_0 + \beta_1 x_{i,1} + \cdots + \beta_K x_{i,K} + \varepsilon_i = \alpha_0 + \mathbf{x}_i^t \boldsymbol{\beta} + \varepsilon_i, \qquad i = 1, \ldots, n.$

Received July 2003; revised April 2004.
[1]Supported by NSF Grant DMS-04-05675.
[2]Supported by NIH Grant K25-CA89867 and NSF Grant DMS-04-05072.
*AMS 2000 subject classifications.* Primary 62J07; secondary 62J05.
*Key words and phrases.* Generalized ridge regression, hypervariance, model averaging, model uncertainty, ordinary least squares, penalization, rescaling, shrinkage, stochastic variable selection, Zcut.

The $\varepsilon_i$ are independent random variables (but not necessarily identically distributed) such that $\mathbb{E}(\varepsilon_i) = 0$ and $\mathbb{E}(\varepsilon_i^2) = \sigma^2$. The variance $\sigma^2 > 0$ is assumed to be unknown.

The true value for $\boldsymbol{\beta}$ will be denoted by $\boldsymbol{\beta}_0 = (\beta_{1,0}, \ldots, \beta_{K,0})^t$ and the true variance of $\varepsilon_i$ by $\sigma_0^2 > 0$. The complexity, or true dimension, is the number of $\beta_{k,0}$ coefficients that are nonzero, which we denote by $k_0$. We assume that $1 \leq k_0 \leq K$, where $K$, the total number of covariates, is a fixed value. For convenience, and without loss of generality, we assume that covariates have been centered and rescaled so that $\sum_{i=1}^{n} x_{i,k} = 0$ and $\sum_{i=1}^{n} x_{i,k}^2 = n$ for each $k = 1, \ldots, K$. Because we can define $\alpha_0 = \bar{Y}$, the mean of the $Y_i$ responses, and replace $Y_i$ by the centered values $Y_i - \bar{Y}$, we can simply assume that $\alpha_0 = 0$. Thus, we remove $\alpha_0$ throughout our discussion.

The classical variable selection framework involves identification of the nonzero elements of $\boldsymbol{\beta}_0$ and sometimes, additionally, estimation of $k_0$. Information-theoretic approaches [see, e.g., Shao (1997)] consider all $2^K$ models and select the model with the best fit according to some information based criteria. These have been shown to have optimal asymptotic properties, but finite sample performance has suffered [Bickel and Zhang (1992), Rao (1999), Shao and Rao (2000) and Leeb and Pötscher (2003)]. Furthermore, such methods become computationally infeasible even for relatively small $K$. Some solutions have been proposed [see, e.g., Zheng and Loh (1995, 1997)] where a data-based ordering of the elements of $\boldsymbol{\beta}$ is used in tandem with a complexity recovery criterion. Unfortunately, the asymptotic rates that need to be satisfied serve only as a guide and can prove difficult to implement in practice.

Bayesian spike and slab approaches to variable selection (see Section 2) have also been proposed [Mitchell and Beauchamp (1988), George and McCulloch (1993), Chipman (1996), Clyde, DeSimone and Parmigiani (1996), Geweke (1996) and Kuo and Mallick (1998)]. These involve designing a hierarchy of priors over the parameter and model spaces of (1). Gibbs sampling is used to identify promising models with high posterior probability of occurrence. The choice of priors is often tricky, although empirical Bayes approaches can be used to deal with this issue [Chipman, George and McCulloch (2001)]. With increasing $K$, however, the task becomes more difficult. Furthermore, Barbieri and Berger (2004) have shown that in many circumstances the high frequency model is not the optimal predictive model and that the median model (the model consisting of those variables which have overall posterior inclusion probability greater than or equal to 50%) is predictively optimal.

In recent work, Ishwaran and Rao (2000, 2003, 2005) used a modified rescaled spike and slab model that makes use of continuous bimodal priors for hypervariance parameters (see Section 3). This method proved particularly suitable for regression settings with very large $K$. Applications of this work included identifying differentially expressing genes from DNA microarray data. It was shown that this could be cast as a special case of (1) under a near orthogonal design for two group problems [Ishwaran and Rao (2003)], and as an orthogonal design for general multiclass problems [Ishwaran and Rao (2005)]. Along the lines of Barbieri and Berger (2004), attention was focused on processing posterior information for $\boldsymbol{\beta}$ (in this case by considering posterior mean values) rather than finding high frequency models. This is because in high-dimensional situations it is common for there to be no high frequency model (in the microarray examples considered $K$ was on the order of 60,000). Improved performance was observed over traditional methods and attributed to the procedure's ability to maintain a balance between low false detection and high statistical power. A partial theoretical analysis was carried out and connections to frequentist shrinkage made. The improved performance was linked to selective shrinkage in which only truly zero coefficients were shrunk toward zero from their ordinary least squares (OLS) estimates. In addition, a novel shrinkage plot which allowed adaptive calibration of significance levels to account for multiple testing under the large $K$ setup was developed.

1.1. *Statement of main results.* In this article we provide a general analysis of the spike and slab approach. A key ingredient to our approach involves drawing upon connections between the posterior mean, the foundation of our variable selection approach, and frequentist generalized ridge regression estimation. Our primary findings are summarized as follows:

1. The use of a spike and slab model with a continuous bimodal prior for hypervariances has distinct advantages in terms of calibration. However, like any prior, its effect becomes swamped by the likelihood as the sample size $n$ increases, thus reducing the potential for the prior to impact model selection relative to a frequentist method. Instead, we introduce a rescaled spike and slab model defined by replacing the $Y$-responses with $\sqrt{n}$-rescaled values. This makes it possible for the prior to have a nonvanishing effect, and so is a type of sample size universality for the prior.
2. This rescaling is accompanied by a variance inflation parameter $\lambda_n$. It is shown through the connection to generalized ridge regression that $\lambda_n$ controls the amount of shrinkage the posterior mean exhibits relative to the OLS, and thus can be viewed as a penalization effect. Theorem 2 of Section 3 shows that if $\lambda_n$ satisfies $\lambda_n \to \infty$ and $\lambda_n/n \to 0$, then the effect of shrinkage vanishes asymptotically and the posterior mean (after suitable rescaling) is asymptotically equivalent to the OLS (and, therefore, is consistent for $\boldsymbol{\beta}_0$).
3. While consistency is important from an estimation perspective, we show for model selection purposes that the most interesting case occurs when $\lambda_n = n$. At this level of penalization, at least for orthogonal designs, the posterior mean achieves an oracle risk misclassification performance relative to the OLS under a correctly chosen value for the hypervariance (Theorem 5 of Section 5). While this is an oracle result, we show that similar risk performance is achieved using a continuous bimodal prior. Continuity of the prior will be shown to be essential for the posterior mean to identify nonzero coefficients, while bimodality of the prior will enable the posterior mean to identify zero coefficients (Theorem 6 of Section 5).

4. Thus, the use of a rescaled spike and slab model, in combination with a continuous bimodal prior, has the effect of turning the posterior mean into a highly effective Bayesian test statistic. Unlike the analogous frequentist test statistic based on the OLS, the posterior mean takes advantage of model averaging and the benefits of shrinkage through generalized ridge regression estimation. This leads to a type of "selective shrinkage" where the posterior mean is asymptotically biased and shrunk toward zero for coefficients that are zero (see Theorem 6 for an explicit finite sample description of the posterior). The exact nature of performance gains compared to standard OLS model selection procedures has to do primarily with this selective shrinkage.

5. Information from the posterior could be used in many ways to select variables; however, by using a local asymptotic argument, we show that the posterior is asymptotically maximized by the posterior mean (see Section 4). This naturally suggests the use of the posterior mean, especially when combined with a reliable thresholding rule. Such a rule, termed "Zcut", is motivated by a ridge distribution that appears in the limit in our analysis. Also suggested from this analysis is a new multivariate null distribution for testing if a coefficient is zero (Section 5).

6. We introduce a forward stepwise selection strategy as an empirical tool for verifying the ability of a model averaging procedure to reduce model uncertainty. If a procedure is effective, then its data based version of the forward stepwise procedure should outperform an OLS model estimator. See Section 6 and Theorem 8.

1.2. *Selective shrinkage.* A common thread underlying the article, and key to most of the results just highlighted, is the selective shrinkage ability of the posterior. It is worthwhile, therefore, to briefly amplify what we mean by this. Figure 1 serves as an illustration of the idea. There $Z$-test statistics $\widehat{\widetilde{Z}}_{k,n}$ estimated by OLS under the full model are plotted against the corresponding posterior mean values $\widehat{\beta}^*_{k,n}$ under our rescaled spike and slab model (the notation used will be explained later in the paper). As mentioned, these rescaled models are derived under a $\sqrt{n}$-rescaling of the data, which forces the posterior mean onto a $\sqrt{n}$-scale. This is why we plot the posterior mean against a test statistic. The results depicted in Figure 1 are based on a simulation, as in Breiman (1992), for an uncorrelated (near-orthogonal) design where $k_0 = 105$, $K = 400$ and $n = 800$ (see Section 8 for details). Selective shrinkage has to do with shrinkage for the zero $\beta_{k,0}$ coefficients, and is immediately obvious from Figure 1. Note how the $\widehat{\beta}^*_{k,n}$ are shrunken toward a cluster of values near zero for many of the zero coefficients, but are similar to the frequentist $Z$-tests for many of the nonzero coefficients. It is precisely this effect we refer to as selective shrinkage.

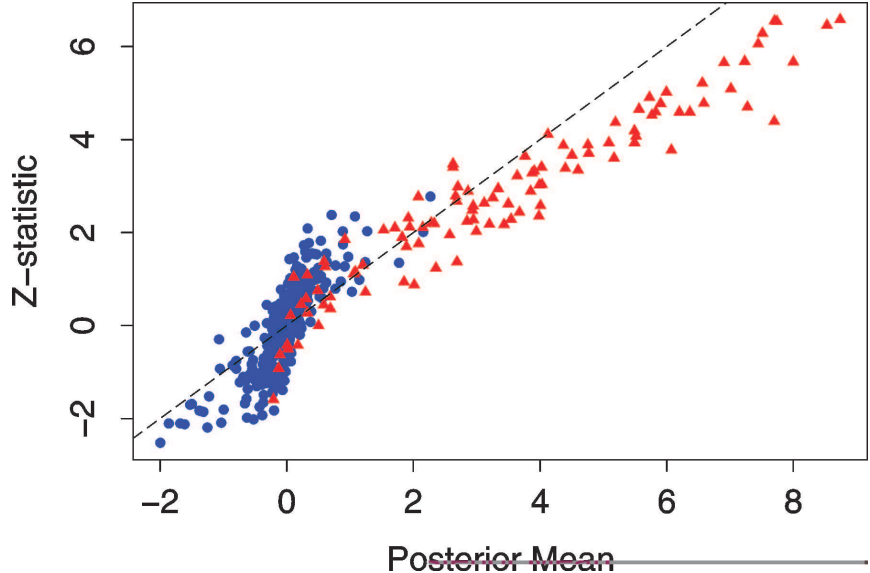

Fig. 1. *Selective shrinkage. Z-test statistics $\widehat{Z}_{k,n}$ versus posterior mean values $\widehat{\beta}^*_{k,n}$ (blue circles are zero coefficients, red triangles nonzero coefficients). Result from Breiman simulation of Section 8 with an uncorrelated design matrix, $k_0 = 105$, $K = 400$ and $n = 800$.*

In fact, this kind of selective shrinkage is not unique to the Bayesian variable selection framework. Shao (1993, 1996) and Zhang (1993) studied cross-validation and bootstrapping for model selection and discovered that to achieve optimal asymptotic performance, a nonvanishing bias term was needed, and this could be constructed by modifying the resampling scheme (see the references for details). Overfit models (ones with too many parameters) are preferentially selected without this bias term. As a connection to this current work, this amounts to detecting zero coefficients—which is a type of selective shrinkage.

1.3. *Organization of the article.* The article is organized as follows. Section 2 presents an overview of spike and slab models. Section 3 introduces our rescaled models and discusses the universality of priors, the role of rescaling and generalized ridge regression. Section 4 examines the optimality of the posterior mean under a local asymptotics framework. Section 5 introduces the Zcut selection strategy. Its optimality in terms of risk performance and complexity recovery is discussed. Section 6 uses a special paradigm in which $\boldsymbol{\beta}_0$ is assumed ordered a priori, and derives both forward and backward selection strategies in the spirit of Leeb and Pötscher (2003). These are used to study the effects of model uncertainty. Sections 7 and 8 present a real data analysis and simulation.

**2. Spike and slab models.** By a *spike and slab model* we mean a Bayesian model specified by the following prior hierarchy:

$$
\begin{aligned}
(Y_i|\mathbf{x}_i, \boldsymbol{\beta}, \sigma^2) &\stackrel{\text{ind}}{\sim} \mathrm{N}(\mathbf{x}_i^t \boldsymbol{\beta}, \sigma^2), \qquad i=1,\ldots,n, \\
(\boldsymbol{\beta}|\boldsymbol{\gamma}) &\sim \mathrm{N}(\mathbf{0}, \boldsymbol{\Gamma}), \\
\boldsymbol{\gamma} &\sim \pi(d\boldsymbol{\gamma}), \\
\sigma^2 &\sim \mu(d\sigma^2),
\end{aligned}
\tag{2}
$$

where $\mathbf{0}$ is the $K$-dimensional zero vector, $\boldsymbol{\Gamma}$ is the $K \times K$ diagonal matrix $\mathrm{diag}(\gamma_1, \ldots, \gamma_K)$, $\pi$ is the prior measure for $\boldsymbol{\gamma} = (\gamma_1, \ldots, \gamma_K)^t$ and $\mu$ is the prior measure for $\sigma^2$. Throughout we assume that both $\pi$ and $\mu$ are chosen to exclude values of zero with probability one; that is, $\pi\{\gamma_k > 0\} = 1$ for $k=1,\ldots,K$ and $\mu\{\sigma^2 > 0\} = 1$.

Lempers (1971) and Mitchell and Beauchamp (1988) were among the earliest to pioneer the spike and slab method. The expression "spike and slab" referred to the prior for $\boldsymbol{\beta}$ used in their hierarchical formulation. This was chosen so that $\beta_k$ were mutually independent with a two-point mixture distribution made up of a uniform flat distribution (the slab) and a degenerate distribution at zero (the spike). Our definition (2) deviates significantly from this. In place of a two-point mixture distribution, we assume that $\boldsymbol{\beta}$ has a multivariate normal scale mixture distribution specified through the prior $\pi$ for the hypervariance $\boldsymbol{\gamma}$. Our basic idea, however, is similar in spirit to the Lempers–Mitchell–Beauchamp approach. To select variables, the idea is to zero out $\beta_k$ coefficients that are truly zero by making their posterior mean values small. The spike and slab hierarchy (2) accomplishes this through the values for the hypervariances. Small hypervariances help to zero out coefficients, while large values inflate coefficients. The latter coefficients are the ones we would like to select in the final model.

EXAMPLE 1 (*Two-component indifference priors*). A popular version of the spike and slab model, introduced by George and McCulloch (1993), identifies zero and nonzero $\beta_k$'s by using zero–one latent variables $\mathscr{I}_k$. This identification is a consequence of the prior used for $\beta_k$, which is assumed to be a scale mixture of two normal distributions:

$$(\beta_k|\mathscr{I}_k) \stackrel{\text{ind}}{\sim} (1-\mathscr{I}_k)\mathrm{N}(0,\tau_k^2) + \mathscr{I}_k\mathrm{N}(0,c_k\tau_k^2), \qquad k=1,\ldots,K.$$

[We use the notation $\mathrm{N}(0,v^2)$ informally here to represent the measure of a normal variable with mean 0 and variance $v^2$.] The value for $\tau_k^2 > 0$ is some suitably small value while $c_k > 0$ is some suitably large value. Coefficients that are promising have posterior latent variables $\mathscr{I}_k = 1$. These coefficients will have large posterior hypervariances and, consequently, large posterior $\beta_k$ values. The opposite occurs when $\mathscr{I}_k = 0$. The prior hierarchy for $\boldsymbol{\beta}$ is completed by assuming a prior for $\mathscr{I}_k$. In principle, one can use any prior over the $2^K$ possible values for $(\mathscr{I}_1,\ldots,\mathscr{I}_K)^t$; however, often $\mathscr{I}_k$ are taken as independent Bernoulli$(w_k)$ random variables, where $0 < w_k < 1$. It is common practice to set $w_k = 1/2$. This is referred to as an *indifference*, or uniform prior. It is clear this setup can be recast as a spike and slab model (2). That is, the prior $\pi(d\boldsymbol{\gamma})$ in (2) is defined by the conditional distributions

(3)
$$(\gamma_k|c_k,\tau_k^2,\mathscr{I}_k) \stackrel{\text{ind}}{\sim} (1-\mathscr{I}_k)\delta_{\tau_k^2}(\cdot) + \mathscr{I}_k \delta_{c_k\tau_k^2}(\cdot), \qquad k=1,\ldots,K,$$
$$(\mathscr{I}_k|w_k) \stackrel{\text{ind}}{\sim} (1-w_k)\delta_0(\cdot) + w_k\delta_1(\cdot),$$

where $\delta_v(\cdot)$ is used to denote a discrete measure concentrated at the value $v$. Of course, (3) can be written more compactly as

$$(\gamma_k|c_k,\tau_k^2,w_k) \stackrel{\text{ind}}{\sim} (1-w_k)\delta_{\tau_k^2}(\cdot) + w_k\delta_{c_k\tau_k^2}(\cdot), \qquad k=1,\ldots,K.$$

However, (3) is often preferred for computational purposes.

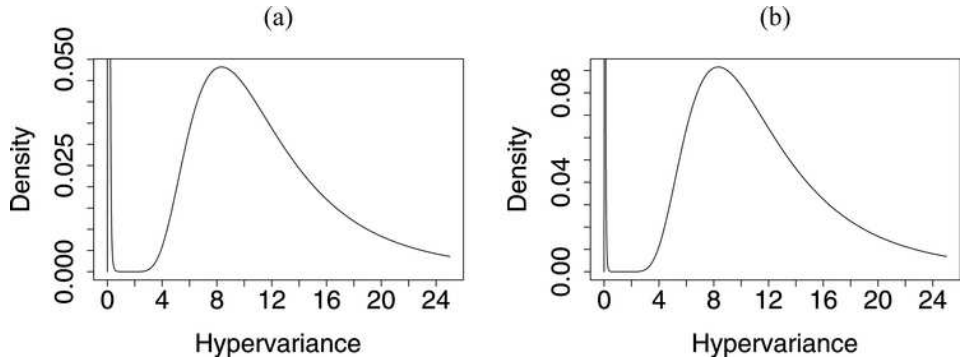

Fig. 2. *Conditional density for $\gamma_k$, where $v_0 = 0.005$, $a_1 = 5$ and $a_2 = 50$ and (a) $w = 0.5$, (b) $w = 0.95$. Observe that only the height of the density changes as $w$ is varied. [Note as $w$ has a uniform prior, (a) also corresponds to the marginal density for $\gamma_k$.]*

EXAMPLE 2 (*Continuous bimodal priors*). In practice, it can be difficult to select the values for $\tau_k^2$, $c_k\tau_k^2$ and $w_k$ used in the priors for $\boldsymbol{\beta}$ and $\mathscr{I}_k$. Improperly chosen values lead to models that concentrate on either too few or too many coefficients. Recognizing this problem, Ishwaran and Rao (2000) proposed a continuous bimodal distribution for $\gamma_k$ in place of the two-point mixture distribution for $\gamma_k$ in (3). They introduced the following prior hierarchy for $\boldsymbol{\beta}$:

$$
\begin{aligned}
(\beta_k | \mathscr{I}_k, \tau_k^2) &\stackrel{\text{ind}}{\sim} \text{N}(0, \mathscr{I}_k \tau_k^2), \qquad k = 1, \ldots, K, \\
(\mathscr{I}_k | v_0, w) &\stackrel{\text{i.i.d.}}{\sim} (1-w)\,\delta_{v_0}(\cdot) + w\,\delta_1(\cdot), \\
(\tau_k^{-2} | a_1, a_2) &\stackrel{\text{i.i.d.}}{\sim} \text{Gamma}(a_1, a_2), \\
w &\sim \text{Uniform}[0, 1].
\end{aligned}
\tag{4}
$$

The prior $\pi$ for $\boldsymbol{\gamma}$ is induced by $\gamma_k = \mathscr{I}_k \tau_k^2$, and thus integrating over $w$ shows that (4) is a prior for $\boldsymbol{\beta}$ as in (2).

In (4), $v_0$ (a small near zero value) and $a_1$ and $a_2$ (the shape and scale parameters for a gamma density) are chosen so that $\gamma_k = \mathscr{I}_k \tau_k^2$ has a continuous bimodal distribution with a spike at $v_0$ and a right continuous tail (see Figure 2). The spike at $v_0$ is important because it enables the posterior to shrink values for the zero $\beta_{k,0}$ coefficients, while the right-continuous tail is used to identify nonzero parameters. Continuity is crucial because it avoids having to manually set a bimodal prior as in (3). Another unique feature of (4) is the parameter $w$. Its value controls how likely $\mathscr{I}_k$ equals 1 or $v_0$, and, therefore, it takes on the role of a *complexity parameter* controlling the size of models. Notice that using an indifference prior is equivalent to choosing a degenerate prior for $w$ at the value of $1/2$. Using a continuous prior for $w$, therefore, allows for a greater amount of adaptiveness in estimating model size.

**3. Rescaling, penalization and universality of priors.** The flexibility of a prior like (4) greatly simplifies the problems of calibration. However, just like any other prior, its effect on the posterior vanishes as the sample size increases, and without some basic adjustment to the underlying spike and slab model, the only way to avoid a washed out effect would be to tune the prior as a function of the sample size. Having to adjust the prior is undesirable. Instead, to achieve a type of "universality," or sample size invariance, we introduce a modified rescaled spike and slab model (Section 3.1). This involves replacing the original $Y_i$ values with ones transformed by a $\sqrt{n}$ factor. Also included in the models is a variance inflation factor needed to adjust to the new variance of the transformed data. To determine an appropriate choice for the inflation factor, we show that this value can also be interpreted as a penalization shrinkage effect of the posterior mean. We show that a value of $n$ is the most appropriate because it ensures that the prior has a nonvanishing effect. This is important, because as we demonstrate later in Section 5, this nonvanishing effect, in combination with an appropriately selected prior for $\boldsymbol{\gamma}$, such as (4), yields a model selection procedure based on the posterior mean with superior performance over one using the OLS.

For our results we require some fairly mild constraints on the behavior of covariates.

*Design assumptions.* Let $\mathbf{X}$ be the $n \times K$ design matrix from the regression model (1). We shall make use of one, or several, of the following conditions:

(D1) $\sum_{i=1}^{n} x_{i,k} = 0$ and $\sum_{i=1}^{n} x_{i,k}^2 = n$ for each $k = 1, \ldots, K$.
(D2) $\max_{1 \leq i \leq n} \|\mathbf{x}_i\|/\sqrt{n} \to 0$, where $\|\cdot\|$ is the $\ell_2$-norm.
(D3) $\mathbf{X}^t\mathbf{X}$ is positive definite.
(D4) $\boldsymbol{\Sigma}_n = \mathbf{X}^t\mathbf{X}/n \to \boldsymbol{\Sigma}_0$, where $\boldsymbol{\Sigma}_0$ is positive definite.

Condition (D1) simply reiterates the assumption that covariates are centered and rescaled. Condition (D2) is designed to keep any covariate $\mathbf{x}_i$ from becoming too large. Condition (D3) will simplify some arguments, but is unnecessary for asymptotic results in light of condition (D4). Condition (D3) is convenient, because it frees us from addressing noninvertibility for small values of $n$. It also allows us to write out closed form expressions for the OLS estimate without having to worry about generalized inverses. Note, however, that from a practical point of view, noninvertibility for $\boldsymbol{\Sigma}_n$ is not problematic. This is because the conditional posterior mean is a generalized ridge estimator, which always exists if the ridge parameters are nonzero.

REMARK 1. We call $\widehat{\boldsymbol{\beta}}_{\mathrm{R}}$ a *generalized ridge estimator* for $\boldsymbol{\beta}_0$ if $\widehat{\boldsymbol{\beta}}_{\mathrm{R}} = (\mathbf{X}^t\mathbf{X} + \mathbf{D})^{-1}\mathbf{X}^t\mathbf{Y}$, where $\mathbf{D}$ is a $K \times K$ diagonal matrix. Here $\mathbf{Y} = (Y_1, \ldots, Y_n)^t$ is the vector of responses. The diagonal elements $d_1, \ldots, d_K$ of $\mathbf{D}$ are assumed to be nonnegative and are referred to as the ridge parameters, while $\mathbf{D}$ is referred to as the ridge matrix. If $d_k > 0$ for each $k$, then $\mathbf{X}^t\mathbf{X} + \mathbf{D}$ is always of full rank. See Hoerl (1962) and Hoerl and Kennard (1970) for background on ridge regression.

3.1. *Rescaled spike and slab models.* By a *rescaled spike and slab model*, we mean a spike and slab model modified as follows:

$$
\begin{aligned}
(Y_i^* | \mathbf{x}_i, \boldsymbol{\beta}, \sigma^2) &\stackrel{\mathrm{ind}}{\sim} \mathrm{N}(\mathbf{x}_i^t \boldsymbol{\beta}, \sigma^2 \lambda_n), \quad i = 1, \ldots, n, \\
(\boldsymbol{\beta} | \boldsymbol{\gamma}) &\sim \mathrm{N}(\mathbf{0}, \boldsymbol{\Gamma}), \\
\boldsymbol{\gamma} &\sim \pi(d\boldsymbol{\gamma}), \\
\sigma^2 &\sim \mu(d\sigma^2),
\end{aligned}
\tag{5}
$$

where $Y_i^* = \widehat{\sigma}_n^{-1} n^{1/2} Y_i$ are rescaled $Y_i$ values, $\widehat{\sigma}_n^2 = \|\mathbf{Y} - \mathbf{X}\widehat{\boldsymbol{\beta}}_n^\circ\|^2/(n-K)$ is the unbiased estimator for $\sigma_0^2$ based on the full model and $\widehat{\boldsymbol{\beta}}_n^\circ = (\mathbf{X}^t\mathbf{X})^{-1}\mathbf{X}^t\mathbf{Y}$ is the OLS estimate for $\boldsymbol{\beta}_0$ from (1).

The parameter $\lambda_n$ appearing in (5) is one of the key differences between (5) and our earlier spike and slab model (2). One way to think about this value is that it's a variance inflation factor introduced to compensate for the scaling of the $Y_i$'s. Given that a $\sqrt{n}$-scaling is used, the most natural choice for $\lambda_n$ would be $n$, reflecting the correct increase in the variance of the data. However, another way to motivate this choice is through a penalization argument. We show that $\lambda_n$ controls the amount of shrinkage and that a value of $\lambda_n = n$ is the amount of penalization required in order to ensure a shrinkage effect in the limit.

REMARK 2. When $\lambda_n = n$, we have found that $\sigma^2$ in (5) plays an important adaptive role in adjusting the penalty $\lambda_n$, but only by some small amount. Our experience has shown that under this setting the posterior for $\sigma^2$ will concentrate around the value of one, thus fine tuning the amount of penalization. Some empirical evidence of this will be provided later on in Section 8.

REMARK 3. Throughout the paper when illustrating the spike and slab methodology, we use the continuous bimodal priors (4) in tandem with the rescaled spike and slab model (5) under a penalization $\lambda_n = n$. Specifically, we use the model

$$
\begin{aligned}
(Y_i^*|\mathbf{x}_i, \boldsymbol{\beta}, \sigma^2) &\stackrel{\text{ind}}{\sim} \mathrm{N}(\mathbf{x}_i^t\boldsymbol{\beta}, \sigma^2 n), \qquad i = 1, \ldots, n, \\
(\beta_k|\mathscr{I}_k, \tau_k^2) &\stackrel{\text{ind}}{\sim} \mathrm{N}(0, \mathscr{I}_k \tau_k^2), \qquad k = 1, \ldots, K, \\
(\mathscr{I}_k|v_0, w) &\stackrel{\text{i.i.d.}}{\sim} (1-w)\delta_{v_0}(\cdot) + w\delta_1(\cdot), \\
(\tau_k^{-2}|b_1, b_2) &\stackrel{\text{i.i.d.}}{\sim} \mathrm{Gamma}(a_1, a_2), \\
w &\sim \mathrm{Uniform}[0,1], \\
\sigma^{-2} &\sim \mathrm{Gamma}(b_1, b_2),
\end{aligned}
$$

(6)

with hyperparameters specified as in Figure 2 and $b_1 = b_2 = 0.0001$. Later theory will show the benefits of using a model like this. In estimating parameters we use the Gibbs sampling algorithm discussed in Ishwaran and Rao (2000). We refer to this method as *Stochastic Variable Selection*, or SVS for short. The SVS algorithm is easily implemented. Because of conjugacy, each of the steps in the Gibbs sampler can be simulated from well-known distributions (see the Appendix for details). In particular, the draw for $\sigma^2$ is from an inverse-gamma distribution, and, in fact, the choice of an inverse-gamma prior for $\sigma^2$ is chosen primarily to exploit this conjugacy. Certainly, however, other priors for $\sigma^2$ could be used. In light of our previous comment, any continuous prior with bounded support should work well as long as the support covers a range of values that includes one. This is important because some of the later theorems (e.g., Theorem 2 of Section 3.3 and Theorem 7 of Section 5.5) require a bounded support for $\sigma^2$. Such assumptions are not unrealistic.

3.2. *Penalization and generalized ridge regression.* To recast $\lambda_n$ as a penalty term, we establish a connection between the posterior mean and generalized ridge regression estimation. This also shows the posterior mean can be viewed as a model averaged shrinkage estimator, providing motivation for its use [see also George (1986) and Clyde, Parmigiani and Vidakovic (1998) for more background and motivation for shrinkage estimators]. Let $\widehat{\boldsymbol{\beta}}_n^*(\boldsymbol{\gamma},\sigma^2) = \mathbb{E}(\boldsymbol{\beta}|\boldsymbol{\gamma},\sigma^2,\mathbf{Y}^*)$ be the conditional posterior mean for $\boldsymbol{\beta}$ from (5). It is easy to verify

$$\widehat{\boldsymbol{\beta}}_n^*(\boldsymbol{\gamma},\sigma^2) = (\mathbf{X}^t\mathbf{X} + \sigma^2\lambda_n\boldsymbol{\Gamma}^{-1})^{-1}\mathbf{X}^t\mathbf{Y}^*$$
$$= \widehat{\sigma}_n^{-1} n^{1/2}(\mathbf{X}^t\mathbf{X} + \sigma^2\lambda_n\boldsymbol{\Gamma}^{-1})^{-1}\mathbf{X}^t\mathbf{Y},$$

where $\mathbf{Y}^* = (Y_1^*,\ldots,Y_n^*)^t$. Thus, $\widehat{\boldsymbol{\beta}}_n^*(\boldsymbol{\gamma},\sigma^2)$ is the ridge solution to a regression of $\mathbf{Y}^*$ on $\mathbf{X}$ with ridge matrix $\sigma^2\lambda_n\boldsymbol{\Gamma}^{-1}$. Let $\widehat{\boldsymbol{\beta}}_n^* = \mathbb{E}(\boldsymbol{\beta}|\mathbf{Y}^*)$ denote the posterior mean for $\boldsymbol{\beta}$ from (5). Then

$$\widehat{\boldsymbol{\beta}}_n^* = \widehat{\sigma}_n^{-1} n^{1/2}\int\{(\mathbf{X}^t\mathbf{X}+\sigma^2\lambda_n\,\boldsymbol{\Gamma}^{-1})^{-1}\mathbf{X}^t\mathbf{Y}\}(\pi\times\mu)(d\boldsymbol{\gamma},d\sigma^2|\mathbf{Y}^*).$$

Hence, $\widehat{\boldsymbol{\beta}}_n^*$ is a weighted average of ridge shrunken estimates, where the adaptive weights are determined from the posteriors of $\boldsymbol{\gamma}$ and $\sigma^2$. In other words, $\widehat{\boldsymbol{\beta}}_n^*$ is an estimator resulting from *shrinkage in combination with model averaging.*

Now we formalize the idea of $\lambda_n$ as a penalty term. Define $\widehat{\boldsymbol{\theta}}_n^*(\boldsymbol{\gamma},\sigma^2) = \widehat{\sigma}_n\widehat{\boldsymbol{\beta}}_n^*(\boldsymbol{\gamma},\sigma^2)/\sqrt{n}$. It is clear $\widehat{\boldsymbol{\theta}}_n^*(\boldsymbol{\gamma},\sigma^2)$ is the ridge solution to a regression of $\mathbf{Y}$ on $\mathbf{X}$ with ridge matrix $\sigma^2\lambda_n\boldsymbol{\Gamma}^{-1}$. A ridge solution can always be recast as an optimization problem, which is a direct way of seeing how $\lambda_n$ plays a penalization role. It is straightforward to show that

$$(7) \qquad \widehat{\boldsymbol{\theta}}_n^*(\boldsymbol{\gamma},\sigma^2) = \arg\min_{\boldsymbol{\beta}}\left\{\|\mathbf{Y}-\mathbf{X}\boldsymbol{\beta}\|^2 + \lambda_n\sum_{k=1}^{K}\sigma^2\gamma_k^{-1}\beta_k^2\right\},$$

which shows clearly that $\lambda_n$ is a penalty term.

REMARK 4. Keep in mind that to achieve this same kind of penalization effect in the standard spike and slab model, (2), requires choosing a prior that depends upon $n$. To see this, note that the conditional posterior mean $\widehat{\boldsymbol{\beta}}_n(\boldsymbol{\gamma},\sigma^2) = \mathbb{E}(\boldsymbol{\beta}|\boldsymbol{\gamma},\sigma^2,\mathbf{Y})$ from (2) is of the form

$$\widehat{\boldsymbol{\beta}}_n(\boldsymbol{\gamma},\sigma^2) = (\mathbf{X}^t\mathbf{X} + \sigma^2\boldsymbol{\Gamma}^{-1})^{-1}\mathbf{X}^t\mathbf{Y}.$$

Multiplying $\widehat{\boldsymbol{\beta}}_n(\boldsymbol{\gamma},\sigma^2)$ by $\sqrt{n}/\widehat{\sigma}_n$ gives $\widehat{\boldsymbol{\beta}}_n^*(\boldsymbol{\gamma},\sigma^2)$, but only if $\sigma^2$ is $O(\lambda_n)$, or if $\boldsymbol{\Gamma}$ has been scaled by $1/\lambda_n$. Either scenario occurs only when the prior depends upon $n$.

3.3. *How much penalization?* The identity (7) has an immediate consequence for the choice of $\lambda_n$, at least from the point of view of estimation. This can be seen by Theorem 1 of Knight and Fu (2000), which establishes consistency for Bridge estimators (ridge estimation being a special case). Their result can be stated in terms of hypervariance vectors with coordinates satisfying $\gamma_1 = \cdots = \gamma_k = \gamma_0$, where $0 < \gamma_0 < \infty$. For ease of notation, we write $\boldsymbol{\gamma} = \gamma_0\mathbf{1}$, where $\mathbf{1}$ is the $K$-dimensional vector with each coordinate equal to one. Theorem 1 of Knight and Fu (2000) implies the following:

THEOREM 1 [Knight and Fu (2000)]. *Suppose that $\varepsilon_i$ are i.i.d. such that $\mathbb{E}(\varepsilon_i) = 0$ and $\mathbb{E}(\varepsilon_i^2) = \sigma_0^2$. If condition* (D4) *holds and $\lambda_n/n \to \lambda_0 \geq 0$, then*

$$\widehat{\boldsymbol{\theta}}_n^*(\gamma_0\mathbf{1},\sigma^2) \xrightarrow{\mathrm{p}} \arg\min_{\boldsymbol{\beta}}\left\{(\boldsymbol{\beta}-\boldsymbol{\beta}_0)^t\boldsymbol{\Sigma}_0(\boldsymbol{\beta}-\boldsymbol{\beta}_0) + \lambda_0\sigma^2\gamma_0^{-1}\sum_{k=1}^K \beta_k^2\right\}.$$

*In particular, if $\lambda_0 = 0$, then $\widehat{\boldsymbol{\theta}}_n^*(\gamma_0\mathbf{1},\sigma^2) \xrightarrow{\mathrm{p}} \boldsymbol{\beta}_0$.*

Knight and Fu's result shows there is a delicate balance between the rate at which $\lambda_n$ increases and consistency for $\boldsymbol{\beta}_0$. Any sequence $\lambda_n$ which increases at a rate of $O(n)$ or faster will yield an inconsistent estimator, while any sequence increasing more slowly than $n$ will lead to a consistent procedure.

The following is an analogue of Knight and Fu's result applied to rescaled spike and slab models. Observe that this result does not require $\varepsilon_i$ to be identically distributed. The boundedness assumptions for $\pi$ and $\mu$ stated in the theorem are for technical reasons. In particular, the assumption that $\sigma^2$ remains bounded cannot be removed. It is required for the penalization effect to be completely determined through the value for $\lambda_n$, analogous to Theorem 1 (however, recall from Remark 3 that this kind of assumption is not unrealistic).

THEOREM 2. *Assume that* (1) *holds where $\varepsilon_i$ are independent such that $\mathbb{E}(\varepsilon_i) = 0$ and $\mathbb{E}(\varepsilon_i^2) = \sigma_0^2$. Let $\widehat{\boldsymbol{\theta}}_n^* = \widehat{\sigma}_n\widehat{\boldsymbol{\beta}}_n^*/\sqrt{n}$. Assume that conditions* (D3) *and* (D4) *hold. Also, suppose there exists some* $\eta_0 > 0$ *such that* $\pi\{\gamma_k \geq \eta_0\} = 1$ *for each* $k = 1, \ldots, K$ *and that* $\mu\{\sigma^2 \leq s_0^2\} = 1$ *for some* $0 < s_0^2 < \infty$. *If* $\lambda_n/n \to 0$, *then* $\widehat{\boldsymbol{\theta}}_n^* = \widehat{\boldsymbol{\beta}}_n^\circ + O_p(\lambda_n/n) \xrightarrow{\text{p}} \boldsymbol{\beta}_0$.

**4. Optimality of the posterior mean.** Theorem 2 shows that a penalization effect satisfying $\lambda_n/n \to 0$ yields a posterior mean (after rescaling) that is asymptotically consistent for $\boldsymbol{\beta}_0$. While consistency is certainly crucial for estimation purposes, it could be quite advantageous in terms of model selection if we have a shrinkage effect that does not vanish asymptotically and a posterior mean that behaves differently from the OLS. This naturally suggests penalizations of the form $\lambda_n = n$.

The following theorem (Theorem 3) is a first step in quantifying these ideas. Not only does it indicate more precisely the asymptotic behavior of the posterior for $\boldsymbol{\beta}$, but it also identifies the role that the normal hierarchy plays in shrinkage. An important conclusion is that the optimal way to process the posterior in a local asymptotics framework is by the posterior mean. We then begin a systematic study of the posterior mean (Section 5) and show how this can be used for effective model selection.

For this result we assume $\lambda_n = n$. Note that because of the rescaling of the $Y_i$'s, the posterior is calibrated to a $\sqrt{n}$-scale, and thus some type of reparameterization is needed if we want to consider the asymptotic behavior of the posterior mean. We will look at the case when the true parameter shrinks to **0** at a $\sqrt{n}$-rate. Think of this as a "local asymptotics case." In some aspects these results complement the work in Section 3 of Knight and Fu [(2000)]. See also Le Cam and Yang [(1990), Chapter 5] for more on local asymptotic arguments.

We assume that the true model is

(8) $$Y_{ni} = \mathbf{x}_i^t \boldsymbol{\beta}_n + \varepsilon_{ni}, \qquad i = 1, \ldots, n,$$

where for each $n$, $\varepsilon_{n1}, \ldots, \varepsilon_{nn}$ are independent random variables. The true parameter is $\boldsymbol{\beta}_n = \boldsymbol{\beta}_0/\sqrt{n}$. Let $Y_{ni}^* = \sqrt{n} Y_{ni}$. To model (8) we use a rescaled spike and slab model of the form

(9) $$\begin{aligned}(Y_{ni}^* | \mathbf{x}_i, \boldsymbol{\beta}) &\sim \mathrm{N}(\mathbf{x}_i^t \boldsymbol{\beta}, n), \qquad i = 1, \ldots, n, \\ (\boldsymbol{\beta} | \boldsymbol{\gamma}) &\sim \nu(d\boldsymbol{\beta} | \boldsymbol{\gamma}), \\ \boldsymbol{\gamma} &\sim \pi(d\boldsymbol{\gamma}),\end{aligned}$$

where $\nu(d\boldsymbol{\beta}|\boldsymbol{\gamma})$ is the prior for $\boldsymbol{\beta}$ given $\boldsymbol{\gamma}$. Write $\nu$ for the prior measure for $\boldsymbol{\beta}$, that is, the prior for $\boldsymbol{\beta}$ marginalized over $\boldsymbol{\gamma}$. Let $\nu_n(\cdot | \mathbf{Y}_n^*)$ denote the posterior measure for $\boldsymbol{\beta}$ given $\mathbf{Y}_n^* = (Y_{n1}^*, \ldots, Y_{nn}^*)^t$. For simplicity, and without loss of generality, the following theorem is based on the assumption that $\sigma_0^2$ is known. There is no loss in generality in making such an assumption, because if $\sigma_0^2$ were unknown, we could always rescale $Y_{ni}$ by $\sqrt{n} Y_{ni}/\widehat{\sigma}_n$ and replace $\boldsymbol{\beta}_n$ with $\sigma_0 \boldsymbol{\beta}_0/\sqrt{n}$ as long as $\widehat{\sigma}_n^2 \xrightarrow{\text{p}} \sigma_0^2$. Therefore, for convenience we assume $\sigma_0^2 = 1$ is known.

THEOREM 3. *Assume that $\nu$ has a density $f$ that is continuous and positive everywhere. Assume that (8) is the true regression model, where $\varepsilon_{ni}$ are independent such that $\mathbb{E}(\varepsilon_{ni}) = 0$, $\mathbb{E}(\varepsilon_{ni}^2) = \sigma_0^2 = 1$ and $\mathbb{E}(\varepsilon_{ni}^4) \le M$ for some $M < \infty$. If (D1)–(D4) hold, then for each $\boldsymbol{\beta}_1 \in \mathbb{R}^K$ and each $C > 0$,*

$$
\begin{aligned}
(10) \quad & \log\!\left(\frac{\nu_n(S(\boldsymbol{\beta}_1, C/\sqrt{n})|\mathbf{Y}_n^*)}{\nu_n(S(\boldsymbol{\beta}_0, C/\sqrt{n})|\mathbf{Y}_n^*)}\right) \\
& \stackrel{d}{\rightsquigarrow} \log\!\left(\frac{f(\boldsymbol{\beta}_1)}{f(\boldsymbol{\beta}_0)}\right) - \tfrac{1}{2}(\boldsymbol{\beta}_1 - \boldsymbol{\beta}_0)^t \boldsymbol{\Sigma}_0 (\boldsymbol{\beta}_1 - \boldsymbol{\beta}_0) + (\boldsymbol{\beta}_1 - \boldsymbol{\beta}_0)^t \mathbf{Z},
\end{aligned}
$$

*where $\mathbf{Z}$ has a $\mathrm{N}(\mathbf{0}, \boldsymbol{\Sigma}_0)$ distribution. Here $S(\boldsymbol{\beta}, C)$ denotes a sphere centered at $\boldsymbol{\beta}$ with radius $C > 0$.*

Theorem 3 quantifies the asymptotic behavior of the posterior and its sensitivity to the choice of prior for $\boldsymbol{\beta}$. Observe that the log-ratio posterior probability on the left-hand side of (10) can be thought of as a random function of $\boldsymbol{\beta}_1$. Call this function $\Psi_n(\boldsymbol{\beta}_1)$. Also, the expression on the right-hand side of (10),

$$
(11) \quad -\tfrac{1}{2}(\boldsymbol{\beta}_1 - \boldsymbol{\beta}_0)^t \boldsymbol{\Sigma}_0 (\boldsymbol{\beta}_1 - \boldsymbol{\beta}_0) + (\boldsymbol{\beta}_1 - \boldsymbol{\beta}_0)^t \mathbf{Z},
$$

is a random concave function of $\boldsymbol{\beta}_1$ with a unique maximum at $\boldsymbol{\beta}_0 + \boldsymbol{\Sigma}_0^{-1} \mathbf{Z}$, a $\mathrm{N}(\boldsymbol{\beta}_0, \boldsymbol{\Sigma}_0^{-1})$ random vector. Consider the limit under an improper prior for $\boldsymbol{\beta}$, where $f(\boldsymbol{\beta}_0) = f(\boldsymbol{\beta}_1)$ for each $\boldsymbol{\beta}_1$. Then $\Psi_n(\boldsymbol{\beta}_1)$ converges in distribution to (11), which as we said has a unique maximum at a $\mathrm{N}(\boldsymbol{\beta}_0, \boldsymbol{\Sigma}_0^{-1})$ vector. This is the same limiting distribution for $\sqrt{n}\,\widehat{\boldsymbol{\beta}}_n^\circ$, the rescaled OLS, under the settings of the theorem. Therefore, under a flat prior the posterior behaves similarly to the distribution for the OLS. This is intuitive, because with a noninformative prior there is no ridge parameter, and, therefore, no penalization effect.

On the other hand, consider what happens when $\boldsymbol{\beta}$ has a $\mathrm{N}(\mathbf{0}, \boldsymbol{\Gamma}_0)$ prior. Now the distributional limit of $\Psi_n(\boldsymbol{\beta}_1)$ is

$$
\tfrac{1}{2}\boldsymbol{\beta}_0^t \boldsymbol{\Gamma}_0^{-1} \boldsymbol{\beta}_0 - \tfrac{1}{2}\boldsymbol{\beta}_1^t \boldsymbol{\Gamma}_0^{-1} \boldsymbol{\beta}_1 - \tfrac{1}{2}(\boldsymbol{\beta}_1 - \boldsymbol{\beta}_0)^t \boldsymbol{\Sigma}_0 (\boldsymbol{\beta}_1 - \boldsymbol{\beta}_0) + (\boldsymbol{\beta}_1 - \boldsymbol{\beta}_0)^t \mathbf{Z}.
$$

As a function of $\boldsymbol{\beta}_1$, this is once again concave. However, now the maximum is

$$
\boldsymbol{\beta}_1 = (\boldsymbol{\Sigma}_0 + \boldsymbol{\Gamma}_0^{-1})^{-1}(\boldsymbol{\Sigma}_0 \boldsymbol{\beta}_0 + \mathbf{Z}),
$$

which is a $\mathrm{N}(\mathbf{V}_0^{-1}\boldsymbol{\beta}_0, \mathbf{V}_0^{-1}\boldsymbol{\Sigma}_0^{-1}\mathbf{V}_0^{-1})$ random vector, where $\mathbf{V}_0 = \mathbf{I} + \boldsymbol{\Sigma}_0^{-1}\boldsymbol{\Gamma}_0^{-1}$. Let $Q(\cdot|\boldsymbol{\gamma}_0)$ represent this limiting normal distribution.

The distribution $Q(\cdot|\boldsymbol{\gamma}_0)$ is quite curious. It appears to be a new type of asymptotic ridge limit. The next theorem identifies it as the limiting distribution for the posterior mean.

THEOREM 4. *Assume that $\boldsymbol{\beta}$ has a $\mathrm{N}(\mathbf{0}, \boldsymbol{\Gamma}_0)$ prior for some fixed $\boldsymbol{\Gamma}_0$. Let $\widehat{\boldsymbol{\beta}}_{nn}^*(\boldsymbol{\gamma}_0) = \mathbb{E}(\boldsymbol{\beta}|\boldsymbol{\gamma}_0, \mathbf{Y}_n^*)$ be the posterior mean from (9), where (8) is the true model. Under the same conditions as Theorem 3, we have $\widehat{\boldsymbol{\beta}}_{nn}^*(\boldsymbol{\gamma}_0) \stackrel{\mathrm{d}}{\leadsto} Q(\cdot|\boldsymbol{\gamma}_0)$.*

Theorem 4 shows the importance of the posterior mean when coefficients shrink to zero. In combination with Theorem 3, it shows that in such settings the correct estimator for asymptotically maximizing the posterior *must* be the posterior mean if a normal prior with a fixed hypervariance is used. Notice that the data does not have to be normal for this result to hold.

**5. The Zcut method, orthogonality and model selection performance.** Theorem 4 motivates the use of the posterior mean in settings where coefficients may all be zero and when the hypervariance is fixed, but how does it perform in general, and what are the implications for variable selection? It turns out that under an appropriately specified prior for $\boldsymbol{\gamma}$, the posterior mean from a rescaled spike and slab model exhibits a type of selective shrinkage property, shrinking in estimates for zero coefficients, while retaining large estimated values for nonzero coefficients. This is a key property of immense potential. By using a *hard shrinkage rule*, that is, a threshold rule for setting coefficients to zero, we can take advantage of selective shrinkage to define an effective method for selecting variables. We analyze the theoretical performance of such a hard shrinkage model estimator termed "Zcut." Our analysis will be confined to orthogonal designs (i.e., $\boldsymbol{\Sigma}_n = \boldsymbol{\Sigma}_0 = \mathbf{I}$) for rescaled spike and slab models under a penalization of $\lambda_n = n$. Under these settings we show Zcut possesses an oracle like risk misclassification property when compared to the OLS. Specifically, we show there is an oracle hypervariance $\boldsymbol{\gamma}_0$ which leads to uniformly better risk performance (Section 5.3) and that this type of risk performance is achieved by using a continuous bimodal prior as specified by (4).

5.1. *Hard shrinkage rules and limiting null distributions.* The Zcut procedure (see Section 5.2 for a formal definition) uses a hard shrinkage rule based on a standard normal distribution. Coefficients are set to zero by comparing their posterior mean values to an appropriate cutoff from a standard normal. This rule can be motivated using Theorem 4. This will also indicate an alternative thresholding rule that is an adaptive function of the true coefficients. For simplicity, assume that $\mu\{\sigma^2 = 1\} = 1$. Under the assumptions outlined above, Theorem 4 implies that $\widehat{\boldsymbol{\beta}}_n^*(\boldsymbol{\gamma})$, the conditional posterior mean from (5), is approximately distributed as $Q_n(\cdot|\boldsymbol{\gamma})$, a $\mathrm{N}(\sqrt{n}\mathbf{D}\boldsymbol{\beta}_0/\sigma_0, \mathbf{D}^t\mathbf{D})$ distribution, where $\mathbf{D}$ is the diagonal matrix $\mathrm{diag}(\gamma_1/(1+\gamma_1), \ldots, \gamma_K/(1+\gamma_K))$ (to apply the theorem in the nonlocal asymptotics case, simply replace $\boldsymbol{\beta}_0$ with $\sqrt{n}\boldsymbol{\beta}_0/\sigma_0$). Consequently, the (unconditional) posterior mean $\widehat{\boldsymbol{\beta}}_n^*$ should be approximately distributed as

$$Q_n^*(\cdot) = \int Q_n(\cdot|\boldsymbol{\gamma})\pi(d\boldsymbol{\gamma}|\mathbf{Y}^*).$$

This would seem to suggest that in testing whether a specific coefficient $\beta_{k,0}$ is zero, and, therefore, deciding whether its coefficient estimate should be shrunk to zero, we should compare its posterior mean value $\widehat{\beta}_{k,n}^*$ to the $k$th marginal of $Q_n^*$ under the null $\beta_{k,0} = 0$. Given the complexity of the posterior distribution for $\boldsymbol{\gamma}$, it is tricky to work out what this distribution is exactly. However, in its place we could use

$$Q_{k,\text{null}}^*(\cdot) = \int \mathrm{N}\left(0, \left(\frac{\gamma_k}{1+\gamma_k}\right)^2\right) \pi(d\gamma_k|\mathbf{Y}^*).$$

Notice that this is only an approximation to the true null distribution because $\pi(d\gamma_k|\mathbf{Y}^*)$ does not specifically take into account the null hypothesis $\beta_{k,0} = 0$. Nevertheless, we argue that $Q_{k,\text{null}}^*$ is a reasonable choice. We will also show that a threshold rule based on $Q_{k,\text{null}}^*$ is not that different from the Zcut rule which uses a $\mathrm{N}(0,1)$ reference distribution.

Both rules can be motivated by analyzing how $\pi(d\gamma_k|\mathbf{Y}^*)$ depends upon the true value for the coefficient. First consider what happens when $\beta_{k,0} \neq 0$ and the null is misspecified. Then the posterior will asymptotically concentrate on large $\gamma_k$ values and $\gamma_k/(1+\gamma_k)$ should be concentrated near one (see Theorem 6 later in this section). Therefore, $Q^*_{k,\text{null}}$ will be approximately N(0,1). Also, when $\beta_{k,0} \neq 0$, the $k$th marginal distribution for $Q_n(\cdot|\boldsymbol{\gamma})$ is dominated by the mean, which in this case equals $\sigma_0^{-1}\sqrt{n}\beta_{k,0}\gamma_k/(1+\gamma_k)$. Therefore, if $\gamma_k$ is large, $\widehat{\beta}^*_{k,n}$ is of order

$$\sigma_0^{-1}\sqrt{n}\beta_{k,0} + O_p(1) = \sigma_0^{-1}\sqrt{n}\beta_{k,0}(1 + O_p(1/\sqrt{n})),$$

which shows that the null is likely to be rejected if $\widehat{\beta}^*_{k,n}$ is compared to a N(0,1) distribution.

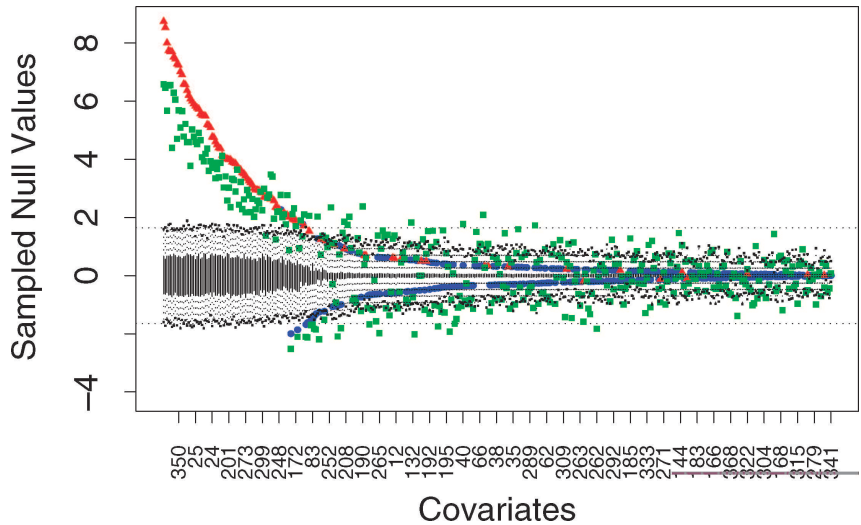

FIG. 3. *Adaptive null intervals. Boxplots of simulated values from $Q_{k,\text{null}}^*$, with $k$ sorted according to largest absolute posterior mean (data from Figure 1). Values from $Q_{k,\text{null}}^*$ were drawn within the SVS Gibbs sampler from a multivariate $\mathrm{N}(\mathbf{0}, \sigma^2 \mathbf{V}_n^t \boldsymbol{\Sigma}_n \mathbf{V}_n)$ distribution where $\mathbf{V}_n = (\boldsymbol{\Sigma}_n + \boldsymbol{\Gamma}^{-1})^{-1}$. Whiskers identify 90% null intervals. Superimposed are $\widehat{Z}_{k,n}$ frequentist test statistics (green squares) and estimated values for $\widehat{\beta}_{k,n}^*$ (blue circles and red triangles used for zero and nonzero coefficients, resp.).*

On the other hand, consider when $\beta_{k,0} = 0$ and the null is really true. Now the hypervariance $\gamma_k$ will often take on small to intermediate values with high posterior probability and $\pi(d\gamma_k|\mathbf{Y}^*)$ should be a good approximation to the posterior under the null. In such settings, using a N(0,1) in place of $Q^*_{k,\text{null}}$ will be slightly more conservative, but this is what we want (after all the null is really true). Let $z_{\alpha/2}$ be the $100 \times (1 - \alpha/2)$ percentile of a standard normal distribution. Observe that

$$\alpha = \mathbb{P}\{|\text{N}(0,1)| \geq z_{\alpha/2}\}$$
$$\geq \int \mathbb{P}\left\{|\text{N}(0,1)| \geq z_{\alpha/2}\left(\frac{\gamma_k}{1+\gamma_k}\right)^{-1}\right\}\pi(d\gamma_k|\mathbf{Y}^*)$$
$$= Q^*_{k,\text{null}}\{|\widehat{\beta}^*_{k,n}| \geq z_{\alpha/2}\}.$$

Therefore, a cut-off value using a N(0,1) distribution yields a significance level larger than $Q^*_{k,\text{null}}$. This is because $Q^*_{k,\text{null}}$ has a smaller variance $\mathbb{E}((\gamma_k/(1+\gamma_k))^2|\mathbf{Y}^*)$ and, therefore, a tighter distribution.

Figure 3 compares the two procedures using the data from our earlier simulation (recall this uses a near orthogonal $\mathbf{X}$ design). Depicted are boxplots for values simulated from $Q^*_{k,\text{null}}$ for each $k$ (see the caption for details). The dashed horizontal lines at $\pm 1.645$ represent a $\alpha = 0.10$ cutoff using a N(0,1) null, while the whiskers for each boxplot are 90% null intervals under $Q^*_{k,\text{null}}$. In this example both procedures lead to similar estimated models, and both yield few false discoveries. In general, however, we prefer the N(0,1) approach because of its simplicity and conservativeness. Nevertheless, the $Q^*_{k,\text{null}}$ intervals can always be produced as part of the analysis. These intervals are valuable because they depict the variability in the posterior mean under the null but are also adaptive to the true value of the coefficient via $\pi(d\gamma_k|\mathbf{Y}^*)$.

5.2. *The Zcut rule.* The preceding argument suggests the use of a thresholding rule that treats the posterior mean as a N(0,1) test statistic. This method, and the resulting hard shrinkage model estimator, have been referred to as Zcut [Ishwaran and Rao (2000, 2003, 2005)]. Here is its formal definition.

THE ZCUT MODEL ESTIMATOR. Let $\widehat{\boldsymbol{\beta}}_n^* = (\widehat{\beta}_{1,n}^*, \ldots, \widehat{\beta}_{K,n}^*)^t$ be the posterior mean for $\boldsymbol{\beta}$ from (5). The Zcut model contains all coefficients $\beta_k$ whose posterior means satisfy $|\widehat{\beta}_{k,n}^*| \geq z_{\alpha/2}$. That is,

$$\text{Zcut} := \{\beta_k : |\widehat{\beta}_{k,n}^*| \geq z_{\alpha/2}\}.$$

Here $\alpha > 0$ is some fixed value specified by the user. The Zcut estimator for $\boldsymbol{\beta}_0$ is the restricted OLS estimator applied to only those coefficients in the Zcut model (all other coefficients are set to zero).

Zcut hard shrinks the posterior mean. Hard shrinkage is important because it reduces the dimension of the model estimator, which is a key to successful subset selection. Given that the posterior mean is already taking advantage of shrinkage, it is natural to wonder how this translates into performance gains over conventional hard shrinkage procedures. We compare Zcut theoretically to "OLS-hard," the hard shrinkage estimator formed from the OLS estimator $\widehat{\boldsymbol{\beta}}_n^\circ = (\widehat{\beta}_{1,n}^\circ, \ldots, \widehat{\beta}_{K,n}^\circ)^t$. Here is its definition:

THE OLS-HARD MODEL ESTIMATOR. The OLS-hard model corresponds to the model with coefficients $\beta_k$ whose Z-statistics, $\widehat{Z}_{k,n}$, satisfy $|\widehat{Z}_{k,n}| \geq z_{\alpha/2}$, where

$$\widehat{Z}_{k,n} = \frac{n^{1/2}\widehat{\beta}_{k,n}^\circ}{\widehat{\sigma}_n(s_{kk})^{1/2}} \tag{12}$$

and $s_{kk}$ is the $k$th diagonal value from $\boldsymbol{\Sigma}_n^{-1}$. That is,

$$\text{OLS-hard} := \{\beta_k : |\widehat{Z}_{k,n}| \geq z_{\alpha/2}\}.$$

The OLS-hard estimator for $\boldsymbol{\beta}_0$ is the restricted OLS estimator using only OLS-hard coefficients.

5.3. *Oracle risk performance.* If Zcut is going to outperform OLS-hard in general, then it is reasonable to expect it will be better in the fixed hypervariance case for some appropriately selected $\boldsymbol{\gamma}$. Theorem 5, our next result, shows this to be true in the context of risk performance. We show there exists a value $\boldsymbol{\gamma} = \boldsymbol{\gamma}_0$ that leads not only to better risk performance, but *uniformly* better risk performance. Let $\mathscr{B}_0 = \{k : \beta_{k,0} = 0\}$ be the indices for the zero coefficients of $\boldsymbol{\beta}_0$. Define

$$\mathscr{R}_Z(\alpha) = \sum_{k \in \mathscr{B}_0} \mathbb{P}\{|\widehat{\beta}_{k,n}^*| \geq z_{\alpha/2}\} + \sum_{k \in \mathscr{B}_0^c} \mathbb{P}\{|\widehat{\beta}_{k,n}^*| < z_{\alpha/2}\}.$$

This is the expected number of coefficients misclassified by Zcut for a fixed $\alpha$-level. This can be thought of as the risk under a zero–one loss function. The misclassification rate for Zcut is $\mathscr{R}_Z(\alpha)/K$. Similarly, define

$$\mathscr{R}_O(\alpha) = \sum_{k \in \mathscr{B}_0} \mathbb{P}\{|\widehat{Z}_{k,n}| \geq z_{\alpha/2}\} + \sum_{k \in \mathscr{B}_0^c} \mathbb{P}\{|\widehat{Z}_{k,n}| < z_{\alpha/2}\}$$

to be the risk for OLS-hard.

THEOREM 5. *Assume that the linear regression model* (1) *holds such that $k_0 < K$ and where $\varepsilon_i$ are i.i.d.* $\mathrm{N}(0, \sigma_0^2)$. *Assume that in* (5) $\boldsymbol{\beta}$ *has a* $\mathrm{N}(\mathbf{0}, \boldsymbol{\Gamma}_0)$ *prior*, $\mu\{\sigma^2 = 1\} = 1$ *and* $\lambda_n = n$. *Then for each $0 < \delta < 1/2$ there exists a $\boldsymbol{\gamma}_0$ such that $\mathscr{R}_Z(\alpha) < \mathscr{R}_O(\alpha)$ for all $\alpha \in [\delta, 1-\delta]$.*

Theorem 5 shows that Zcut's risk is uniformly better than the OLS-hard in *any* finite sample setting if $\boldsymbol{\gamma}$ is set at the oracle value $\boldsymbol{\gamma}_0$. Of course, in practice, this oracle value is unknown, which raises the interesting question of whether the same risk behavior can be achieved by relying on a well-chosen prior for $\boldsymbol{\gamma}$. Also, Theorem 5 requires that $\varepsilon_i$ are normally distributed, but can this assumption be removed?

5.4. *Risk performance for continuous bimodal priors.* Another way to frame these questions is in terms of the posterior behavior of the hypervariances $\gamma_k$. This is because risk performance ultimately boils down to their behavior. One can see this by carefully inspecting the proof of Theorem 5. There the oracle $\boldsymbol{\gamma}_0$ is chosen so that its values are large for the nonzero $\beta_{k,0}$ coefficients and small otherwise. Under any prior $\pi$,

$$\widehat{\beta}_{k,n}^* = \mathbb{E}_\pi\left(\frac{\gamma_k}{1+\gamma_k}\Big|\mathbf{Y}^*\right)\widehat{Z}_{k,n}.$$

In particular, for the $\pi$ obtained by fixing $\boldsymbol{\gamma}$ at $\boldsymbol{\gamma}_0$, the posterior mean is shrunk toward zero for the zero coefficients, thus greatly reducing the number of misclassifications from this group of variables relative to OLS-hard. Meanwhile for the nonzero coefficients, $\widehat{\beta}_{k,n}^*$ is approximately equal to $\widehat{Z}_{k,n}$, so the risk from this group of variables is the same for both procedures, and, therefore, Zcut's risk is smaller overall. Notice that choosing $\boldsymbol{\gamma}_0$ in this fashion also leads to what we have been calling selective shrinkage. So good risk performance follows from selective shrinkage, which ultimately is a statement about the posterior behavior of $\gamma_k$. This motivates the following theorem.

THEOREM 6. *Assume in* (1) *that condition* (D2) *holds and $\varepsilon_i$ are independent such that* $\mathbb{E}(\varepsilon_i) = 0$, $\mathbb{E}(\varepsilon_i^2) = \sigma_0^2$ *and* $\mathbb{E}(\varepsilon_i^4) \leq M$ *for some* $M < \infty$. *Suppose in* (5) *that* $\mu\{\sigma^2 = 1\} = 1$ *and* $\lambda_n = n$.

(a) *If the support for $\pi$ contains a set $[\eta_0, \infty)^K$ for some finite constant $\eta_0 \geq 0$, then, for each small $\delta > 0$,*

$$\pi_n\left(\left\{\boldsymbol{\gamma}: \frac{\gamma_k}{1+\gamma_k} > 1-\delta\right\}\Big|\mathbf{Y}^*\right) \xrightarrow{\mathrm{p}} 1 \qquad \text{if } \beta_{k,0} \neq 0,$$

*where $\pi_n(\cdot|\mathbf{Y}^*)$ is the posterior measure for $\boldsymbol{\gamma}$.*

(b) *Let $f_k^*(\cdot|w)$ denote the posterior density for $\gamma_k$ given $w$. If $\pi$ is the continuous bimodal prior specified by* (4), *then*

$$(13) \quad f_k^*(u|w) \propto \exp\left(\frac{u}{2(1+u)}\xi_{k,n}^2\right)(1+u)^{-1/2}((1-w)g_0(u)+wg_1(u)),$$

*where $g_0(u) = v_0 u^{-2} g(v_0 u^{-1})$, $g_1(u) = u^{-2} g(u^{-1})$,*

$$g(u) = \frac{a_2^{a_1}}{(a_1-1)!} u^{a_1-1} \exp(-a_2 u),$$

*and $\xi_{k,n} = \hat{\sigma}_n^{-1} n^{-1/2} \sum_{i=1}^n x_{i,k} Y_i$. Note that if $\beta_{k,0} = 0$, then $\xi_{k,n} \xrightarrow{\mathrm{d}} \mathrm{N}(0,1)$.*

Part (a) of Theorem 6 shows why continuity for $\pi$ is needed for good risk performance. To be able to selectively shrink coefficients, the posterior must concentrate on arbitrarily large values for the hypervariance when the coefficient is truly nonzero. Part (a) shows this holds asymptotically as long as $\pi$ has an appropriate support. A continuous prior meets this requirement. Selective shrinkage also requires small hypervariances for the zero coefficients, which is what part (b) asserts happens with a continuous bimodal prior. Note importantly that this is a finite sample result and is distribution free. The expression (13) shows that the posterior density for $\gamma_k$ (conditional on $w$) is bimodal. Indeed, except for the leading term

$$(14) \qquad\qquad \exp\left(\frac{u}{2(1+u)}\xi_{k,n}^2\right),$$

which reflects the effect on the prior due to the data, the posterior density is nearly identical to the prior. What (14) does is to adjust the amount of probability at the slab in the prior (cf. Figure 2) using the value of $\xi_{k,n}^2$. As indicated in part (b), if the coefficient is truly zero, then $\xi_{k,n}^2$ will have an approximate $\chi^2$-distribution, so this should introduce a relatively small adjustment. Notice this also implies that the effect of the prior does *not* vanish asymptotically for zero coefficients. This is a key aspect of using a rescaled spike and slab model. Morever, because the posterior for $\gamma_k$ will be similar to the prior when $\beta_{k,0} = 0$, it will concentrate near zero, and hence the posterior mean will be biased and shrunken toward zero relative to the frequentist $Z$-test.

On the other hand, if the coefficient is nonzero, then (14) becomes exponentially large and most of the mass of the density shifts to larger hypervariances. This, of course, matches up with part (a) of the theorem. Figure 4 shows how the posterior cumulative distribution function varies in terms of $\xi_{k,n}^2$. Even for fairly large values of $\xi_{k,n}^2$ (e.g., from the 75th percentile of a $\chi^2$-distribution), the distribution function converges to one rapidly for small hypervariances. This shows that the posterior will concentrate on small hypervariances unless $\xi_{k,n}^2$ is abnormally large.

Figure 5 shows how the hypervariances might vary in a real example. We have plotted the posterior means $\widehat{\beta}_{k,n}^*$ for the Breiman simulation of Figure 1 against $\mathbb{E}((\gamma_k/(1+\gamma_k))^2|\mathbf{Y}^*)$ (the variance of $Q_{k,\text{null}}^*$). This shows quite clearly the posterior's ability to adaptively estimate the hypervariances for selective shrinkage. Figure 6 shows how this selective shrinkage capability is translated into risk performance. As seen, Zcut's misclassification performance is uniformly better than OLS-hard over a wide range of cut-off values, exactly as our theory suggests.

REMARK 5. The assumption in Theorems 5 and 6 that $\mu\{\sigma^2 = 1\} = 1$ is not typical in practice. As discussed, it is beneficial to assume that $\sigma^2$ has a continuous prior to allow adaptive penalization. Nevertheless, Theorem 5 shows even if $\sigma^2 = 1$, thus forgoing the extra benefits of finite sample adaptation, the total risk for Zcut with an appropriately fixed $\boldsymbol{\gamma}_0$ is still uniformly better than OLS-hard. The same argument could be made for Theorem 6. That is, from a theoretical point of view, it is not restrictive to assume a fixed $\sigma^2$.

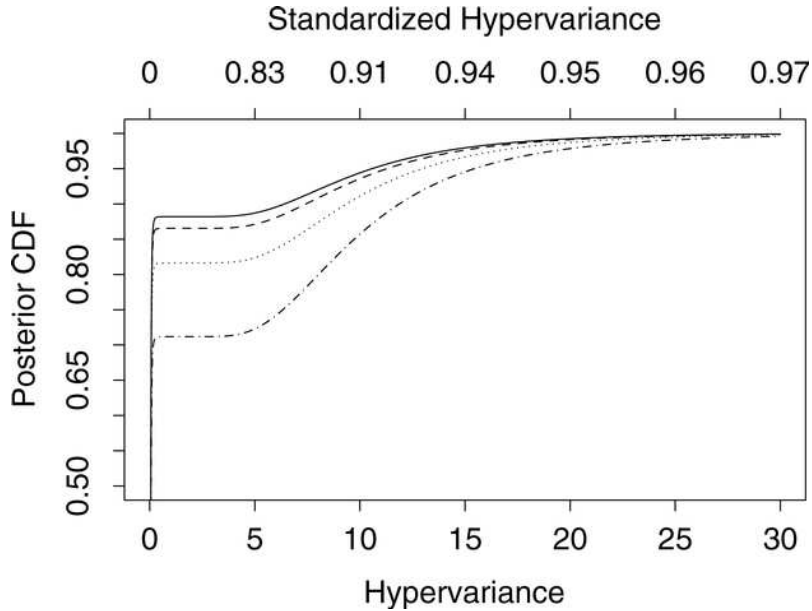

FIG. 4. *Posterior cumulative distribution function for $\gamma_k$ conditional on $w$ (hyperparameters equal to those in Figure 2 and $w = 0.3$). Curves from top to bottom are derived by setting $\xi_{k,n}^2$ at the $25, 50, 75$ and $90$th percentiles for a $\chi^2$-distribution with one degree of freedom. Standardized hypervariance axis defined as $\gamma_k/(1+\gamma_k)$.*

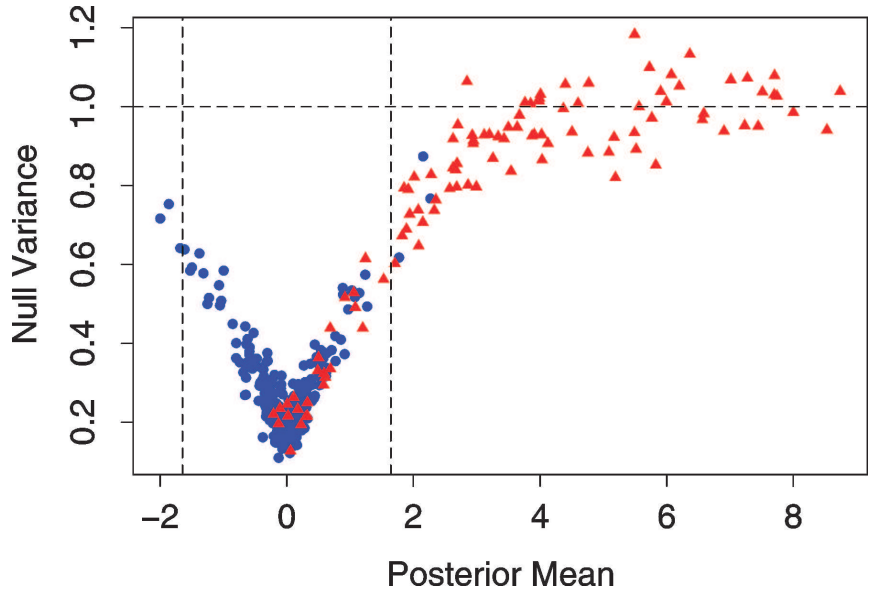

Fig. 5. *Posterior means $\widehat{\beta}^*_{k,n}$ versus variances $\mathbb{E}((\gamma_k/(1+\gamma_k))^2|\mathbf{Y}^*)$ of $Q^*_{k,\text{null}}$ from simulation used in Figure* 1. *Triangles in red are nonzero coefficients.*

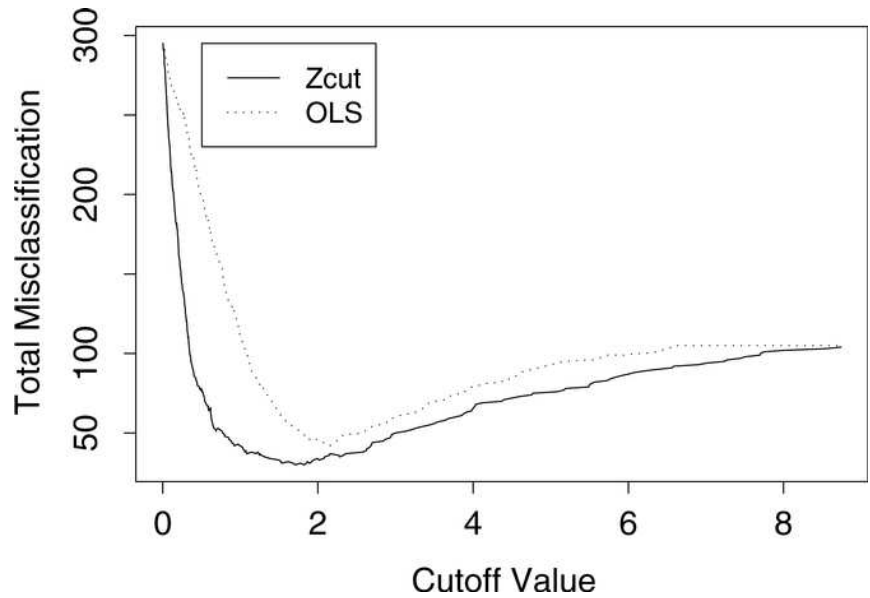

Fig. 6. *Total number of misclassified coefficients from simulation used in Figure* 1. *Observe how Zcut's total misclassification is less than OLS-hard's over a range of cutoff values* $z_{\alpha/2}$.

5.5. *Complexity recovery.* We further motivate Zcut by showing that it consistently estimates the true model under a threshold value that is allowed to change with $n$. Let

$$\mathscr{M}_0 = (\mathbb{I}\{\beta_{1,0} \neq 0\}, \ldots, \mathbb{I}\{\beta_{K,0} \neq 0\})^t$$

be the $K$-dimensional binary vector recording which coordinates of $\boldsymbol{\beta}_0$ are nonzero [$\mathbb{I}(\cdot)$ denotes the indicator function]. By consistent estimation of the true model, we mean the existence of an estimator $\widehat{\mathscr{M}}_n$ such that $\widehat{\mathscr{M}}_n \xrightarrow{\mathrm{P}} \mathscr{M}_0$. We show that such an estimator can be constructed from $\widehat{\boldsymbol{\beta}}_n^*$. Let

$$\mathscr{M}_n(C) = (\mathbb{I}\{|\widehat{\beta}_{1,n}^*| \geq C\}, \ldots, \mathbb{I}\{|\widehat{\beta}_{K,n}^*| \geq C\})^t.$$

The Zcut estimator corresponds to setting $C = z_{\alpha/2}$. The next theorem shows we can consistently recover $\mathscr{M}_0$ by letting $C$ converge to $\infty$ at any rate slower than $\sqrt{n}$.

THEOREM 7. *Assume that the priors $\pi$ and $\mu$ in* (5) *are chosen so that $\pi\{\gamma_k \geq \eta_0\} = 1$ for some $\eta_0 > 0$ for each $k = 1, \ldots, K$ and that $\mu\{\sigma^2 \leq s_0^2\} = 1$ for some $0 < s_0^2 < \infty$. Let $\widehat{\mathscr{M}}_n = \mathscr{M}_n(C_n)$, where $C_n \to \infty$ is any positive increasing sequence such that $C_n/\sqrt{n} \to 0$. Assume that the linear regression model* (1) *holds where $\varepsilon_i$ are independent such that $\mathbb{E}(\varepsilon_i) = 0$, $\mathbb{E}(\varepsilon_i^2) = \sigma_0^2$ and $\mathbb{E}(\varepsilon_i^4) \leq M$ for some $M < \infty$. If* (D2) *holds and $\lambda_n = n$, then $\widehat{\mathscr{M}}_n \xrightarrow{\mathrm{P}} \mathscr{M}_0$.*

An immediate consequence of Theorem 7 is that the true model complexity $k_0$ can be estimated consistently. By the continuous mapping theorem, we obtain the following:

COROLLARY 1. *Let $\widehat{\mathscr{M}}_n = (\widehat{\mathscr{M}}_{1,n}, \ldots, \widehat{\mathscr{M}}_{K,n})^t$ and let $\hat{k}_n = \sum_{k=1}^K \mathbb{I}\{\widehat{\mathscr{M}}_{k,n} \neq 0\}$ be the number of nonzero coordinates of $\widehat{\mathscr{M}}_n$. Then, under the conditions of Theorem 7, $\hat{k}_n \xrightarrow{\mathrm{P}} k_0$.*

**6. The effects of model uncertainty.** In this section we prove an asymptotic complexity result for a specialized type of forward stepwise model selection procedure. This forward stepwise method is a modification of a backward stepwise procedure introduced by Pötscher (1991) and discussed recently in Leeb and Pötscher (2003). We show in orthogonal settings that if the coordinates of $\boldsymbol{\beta}_0$ are perfectly ordered a priori, then the forward stepwise procedure leads to improved complexity recovery relative to the OLS-hard. Interestingly, the backward stepwise procedure has the worst performance of all three methods (Theorem 8 of Section 6.3). This result can be used as an empirical tool for assessing a procedure's ability to reduce model uncertainty. If a model selection procedure is effectively reducing model uncertainty, then it should produce an accurate ranking of coefficients in finite samples. Consequently, the forward stepwise procedure based on this data based ranking should perform better than OLS-hard. This provides an indirect way to confirm a procedure's ability to reduce model uncertainty.

REMARK 6. The idea of pre-ranking covariates and then selecting models has become a well established technique in the literature. As mentioned, this idea was used by Pötscher (1991) and Leeb and Pötscher (2003), but also appears in Zhang (1992), Zheng and Lo (1995, 1997), Rao and Wu (1989) and Ishwaran (2004).

We use this strategy to assess the performance of a rescaled spike and slab model. For a data based ordering of $\boldsymbol{\beta}$, we use the absolute posterior means $|\widehat{\beta}^*_{k,n}|$. The first coordinate of $\boldsymbol{\beta}$ corresponds to the largest absolute posterior value, the second coordinate to the second largest value, and so forth. The data based forward stepwise procedure using this ranking is termed "svsForwd." Section 6.2 provides a formal description. In Section 8 we use simulations to systematically compare the performance of svsForwd to OLS-hard as an indirect way to confirm SVS's ability to reduce model uncertainty. Figure 7 provides some preliminary evidence of this capability. There we have compared a ranking of $\boldsymbol{\beta}$ using the posterior mean against an OLS ordering using $|\widehat{Z}_{k,n}|$. Figure 7 is based on the simulation presented in Figure 1.

We note that it is possible to consistently estimate the order of the $\boldsymbol{\beta}_0$ coordinates using the posterior mean. Let $U_{k,n}$ be the $k$th largest value from the set $\{|\widehat{\beta}^*_{k,n}| : k = 1, \ldots, K\}$. That is, $U_{1,n} \geq U_{2,n} \geq \cdots \geq U_{K,n}$. Let

$$\widehat{\mathscr{M}}_{(n)} = (\mathbb{I}\{U_{1,n} \geq C_n\}, \ldots, \mathbb{I}\{U_{K,n} \geq C_n\})^t,$$

where $C_n$ is a positive sequence satisfying $C_n \to \infty$ and $C_n/\sqrt{n} \to 0$. By inspection of the proof of Theorem 7, Corollary 2 can be shown.

COROLLARY 2. *Under the conditions of Theorem 7, $\widehat{\mathscr{M}}_{(n)} \xrightarrow{\mathrm{p}} (1, \ldots, 1, \mathbf{0}^t_{K-k_0})^t$.*

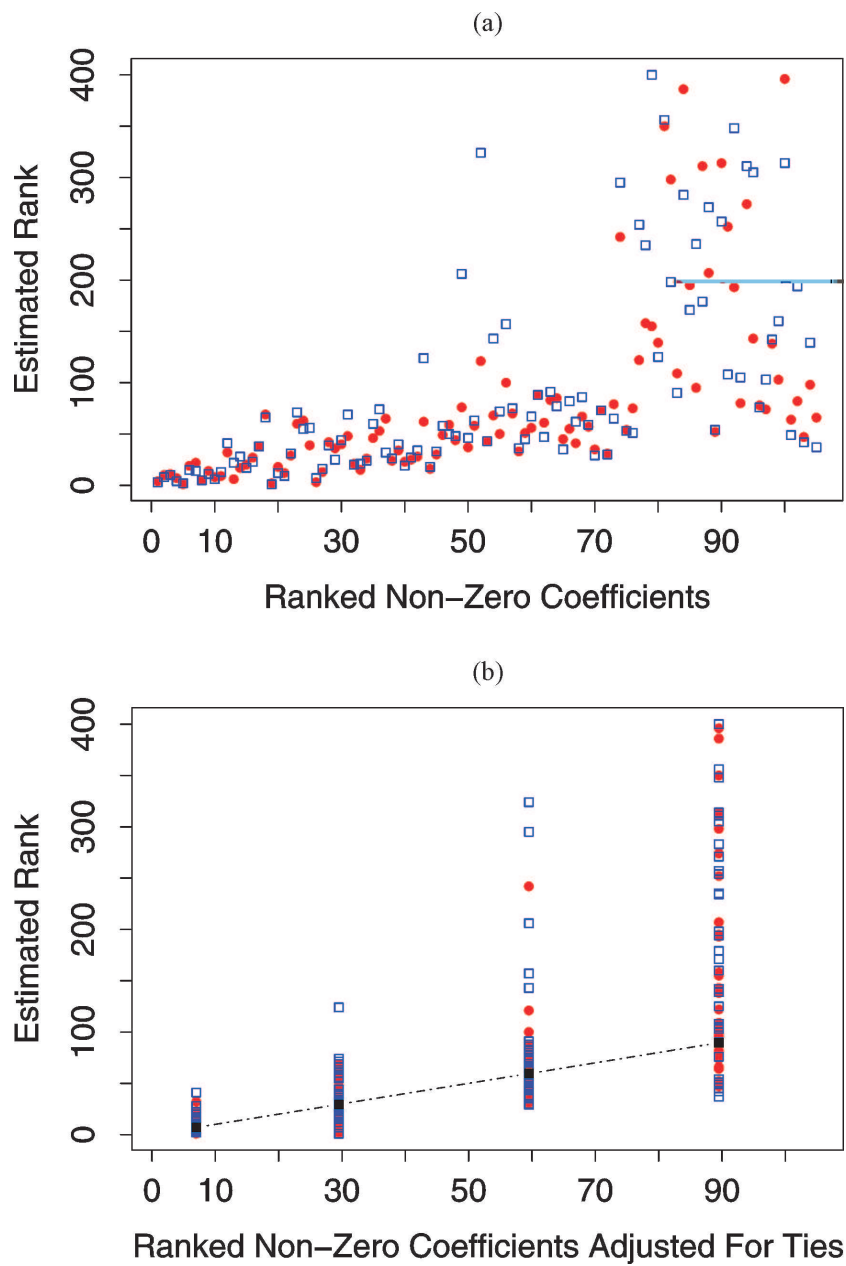

FIG. 7. (a) *True rank of a coefficient versus estimated rank using posterior means (circles) and OLS (squares). The lower the rank, the larger the absolute value of the coefficient. Data from Breiman simulation of Figure* 1 *(only nonzero coefficients shown).* (b) *Same plot as* (a) *but with true ranks averaged to adjust for ties in true coefficient values (simulation used four unique nonzero coefficient values). Dashed line connects values for true average rank. Note the higher variability in OLS, especially for intermediate coefficients.*

6.1. *Backward model selection.* We begin by reviewing the backward stepwise procedure of Pötscher (1991). For notational ease, we avoid subscripts of $n$ as much as possible. We assume the coordinates of $\boldsymbol{\beta}_0$ have been ordered, so that the first $k_0$ coordinates are nonzero. That is,

$$\boldsymbol{\beta}_0 = (\beta_{1,0}, \ldots, \beta_{k_0,0}, \mathbf{0}_{K-k_0}^t)^t,$$

where $\mathbf{0}_{K-k_0}$ is the $(K-k_0)$-dimensional zero vector. We assume the design matrix $\mathbf{X}$ has been suitably recoded as well. Let $\mathbf{X}[k]$ be the $n \times k$ design matrix formed from the first $k$ columns of the re-ordered $\mathbf{X}$. Let

$$\widehat{\boldsymbol{\beta}}^{\circ}[k] = (\widehat{\beta}_1^{\circ}[k], \ldots, \widehat{\beta}_k^{\circ}[k])^t$$
$$= (\mathbf{X}[k]^t \mathbf{X}[k])^{-1} \mathbf{X}[k]^t \mathbf{Y}$$

be the restricted OLS estimator using only the first $k$ variables. To test whether the $k$th coefficient $\beta_k$ is zero, define the test statistic

$$(15) \qquad \widetilde{Z}_{k,n} = \frac{n^{1/2} \widehat{\beta}_k^{\circ}[k]}{\widehat{\sigma}_n (s_{kk}[k])^{1/2}},$$

where $s_{kk}[k]$ is the $k$th diagonal value from $(\mathbf{X}[k]^t \mathbf{X}[k]/n)^{-1}$. Let $\alpha_1, \ldots, \alpha_K$ be a sequence of fixed positive $\alpha$-significance values for the $\widetilde{Z}_{k,n}$ test statistics. Estimate the true complexity $k_0$ by the estimator $\hat{k}_B$, where

$$\hat{k}_B = \max\{k : |\widetilde{Z}_{k,n}| \geq z_{\alpha_k/2}, k = 0, \ldots, K\}.$$

To ensure that $\hat{k}_B$ is well defined, take $\widetilde{Z}_{0,n} = 0$ and $z_{\alpha_0/2} = 0$.

Observe if $\hat{k}_B = k$, then $\widetilde{Z}_{k,n}$ is the first test starting from $k = K$ and going backward to $k = 0$ satisfying $|\widetilde{Z}_{k,n}| \geq z_{\alpha_k/2}$ and $|\widetilde{Z}_{j,n}| < z_{\alpha_j/2}$ for $j = k+1, \ldots, K$. This corresponds to accepting the event $\{\boldsymbol{\beta} : \beta_{k+1} = 0, \ldots, \beta_K = 0\}$, but rejecting $\{\boldsymbol{\beta} : \beta_k = 0, \ldots, \beta_K = 0\}$. The post-model selection estimator for $\boldsymbol{\beta}$ is defined as

$$\widehat{\boldsymbol{\beta}}_B = \mathbf{0}_K \mathbb{I}\{\hat{k}_B = 0\} + \sum_{k=1}^{K} (\widehat{\boldsymbol{\beta}}^{\circ}[k]^t, \mathbf{0}_{K-k}^t)^t \mathbb{I}\{\hat{k}_B = k\}.$$

It should be clear that the estimators $\hat{k}_B$ and $\widehat{\boldsymbol{\beta}}_B$ are derived from a backward stepwise mechanism.

REMARK 7. Observe that $\widetilde{Z}_{k,n}$ uses $\widehat{\sigma}_n^2$, the estimate for $\sigma_0^2$ based on the full model, rather than an estimate based on the first $k$ variables, and so, in this way, is different from a conventional stepwise procedure. The latter estimates are only unbiased if $k \geq k_0$ and can perform quite badly otherwise.

REMARK 8. At first glance it seems the backward procedure requires $K$ regression analyses to compute $\widehat{\boldsymbol{\beta}}^\circ[k]$ for each $k$. This would be expensive for large $K$, requiring a computational effort of $O(\sum_{k=1}^{K} k^3)$. In fact, the whole procedure can be reduced to the problem of finding an orthogonal decomposition of the $\mathbf{X}$ matrix, an $O(K^3)$ operation. This idea rests on the following observations implicit in Lemma A.1 of Leeb and Pötscher (2003). Let

$$\mathbf{P}_k^\perp = \mathbf{I} - \mathbf{X}[k](\mathbf{X}[k]^t \mathbf{X}[k])^{-1} \mathbf{X}[k]^t$$

be the projection onto the orthogonal complement of the linear space spanned by $\mathbf{X}[k]$. Let $\mathbf{x}_{(k)}$ denote the $k$th column vector of $\mathbf{X}$ (thus $\mathbf{X}[k] = [\mathbf{x}_{(1)}, \ldots, \mathbf{x}_{(k)}]$). Define

$$\mathbf{u}_1 = \mathbf{x}_{(1)} \quad \text{and} \quad \mathbf{u}_k = \mathbf{P}_{k-1}^\perp \mathbf{x}_{(k)} \qquad \text{for } k = 2, \ldots, K.$$

One can show that

$$\widehat{\beta}_k^\circ[k] = (\mathbf{u}_k^t \mathbf{u}_k)^{-1} \mathbf{u}_k^t \mathbf{Y}, \qquad k = 1, \ldots, K.$$

Consequently, the backward procedure is equivalent to finding an orthogonal decomposition of $\mathbf{X}$. (Note that this argument shows $\widehat{\beta}_1^\circ[1], \ldots, \widehat{\beta}_K^\circ[K]$ are mutually uncorrelated if $\varepsilon_i$ are independent, $\mathbb{E}(\varepsilon_i) = 0$ and $\mathbb{E}(\varepsilon_i) = \sigma_0^2$. See Lemma A.1 of Leeb and Pötscher (2003). This will be important in the proof of Theorem 8.)

6.2. *Forward model selection.* A forward stepwise procedure and its associated post-model selection estimator for $\boldsymbol{\beta}_0$ can be defined in an analogous way. Define

$$(16) \qquad \hat{k}_F = \min\{k - 1 : |\widetilde{Z}_{k,n}| < z_{\alpha_k/2}, k = 1, \ldots, K+1\},$$

where $\widetilde{Z}_{K+1,n} = 0$ and $\alpha_{K+1} = 0$ are chosen to ensure a well-defined procedure. Observe if $\hat{k}_F = k - 1$, then $\widetilde{Z}_{k,n}$ is the first test statistic such that $|\widetilde{Z}_{k,n}| < z_{\alpha_k/2}$, while $|\widetilde{Z}_{j,n}| \geq z_{\alpha_j/2}$ for $j = 1, \ldots, k - 1$. This corresponds to accepting the event $\{\boldsymbol{\beta} : \beta_1 \neq 0, \ldots, \beta_{k-1} \neq 0\}$, but rejecting $\{\boldsymbol{\beta} : \beta_1 \neq 0, \ldots, \beta_k \neq 0\}$. Note that $\hat{k}_F = 0$ if $|\widetilde{Z}_{1,n}| < z_{\alpha_1/2}$. The post-model selection estimator for $\boldsymbol{\beta}_0$ is

$$(17) \qquad \widehat{\boldsymbol{\beta}}_F = \mathbf{0}_K \mathbb{I}\{\hat{k}_F = 0\} + \sum_{k=1}^{K} (\widehat{\boldsymbol{\beta}}^\circ[k]^t, \mathbf{0}_{K-k}^t)^t \mathbb{I}\{\hat{k}_F = k\}.$$

The data based version of forward stepwise, svsForwd, mentioned earlier is defined as follows:

THE SVSFORWD MODEL ESTIMATOR. Re-order the coordinates of $\boldsymbol{\beta}$ (and the columns of the design matrix $\mathbf{X}$) using the absolute posterior means $\widehat{\beta}^*_{k,n}$ from (5). If $\hat{k}_F \geq 1$, the svsForwd model is defined as

$$\text{svsForwd} := \{\beta_k : k = 1, \ldots, \hat{k}_F\};$$

otherwise, if $\hat{k}_F = 0$, let svsForwd be the null model. Define the svsForwd post-model selection estimator for $\boldsymbol{\beta}_0$ as in (17).

6.3. *Complexity recovery.* The following theorem identifies the limiting distribution for $\hat{k}_B$ and $\hat{k}_F$. It also considers OLS-hard. Let $\hat{k}_O$ denote the OLS-hard complexity estimator (i.e., $\hat{k}_O$ equals the number of parameters in OLS-hard). Part (a) of the following theorem is related to Lemma 4 of Pötscher (1991).

THEOREM 8. *Assume that* (D1)–(D4) *hold for* (1), *where* $\varepsilon_i$ *are independent such that* $\mathbb{E}(\varepsilon_i) = 0$, $\mathbb{E}(\varepsilon_i^2) = \sigma_0^2$ *and* $\mathbb{E}(\varepsilon_i^4) \leq M$ *for some* $M < \infty$. *Let* $k_B$, $k_F$ *and* $k_O$ *denote the limits for* $\hat{k}_B$, $\hat{k}_F$ *and* $\hat{k}_O$, *respectively, as* $n \to \infty$. *For* $1 \leq k \leq K$,

(a) $\mathbb{P}\{k_B = k\} = 0 \times \mathbb{I}\{k < k_0\} + (1 - \alpha_{k_0+1}) \cdots (1 - \alpha_K)\mathbb{I}\{k = k_0\}$
$\quad + \alpha_k(1 - \alpha_{k+1}) \cdots (1 - \alpha_K)\mathbb{I}\{k > k_0\}.$

*Moreover, when* $X$ *has an orthogonal design (i.e.,* $\boldsymbol{\Sigma}_n = \boldsymbol{\Sigma}_0 = \mathbf{I}$),

(b) $\mathbb{P}\{k_F = k\} = 0 \times \mathbb{I}\{k < k_0\} + (1 - \alpha_{k_0+1})\mathbb{I}\{k = k_0\}$
$\quad + (1 - \alpha_{k+1})\alpha_{k_0+1} \cdots \alpha_k \mathbb{I}\{k > k_0\}.$
(c) $\mathbb{P}\{k_O = k\} = 0 \times \mathbb{I}\{k < k_0\}$
$\quad + \mathbb{P}\{B_{k_0+1} + \cdots + B_K = k - k_0\}\mathbb{I}\{k \geq k_0\},$

*where* $\alpha_{K+1} = 0$ *in* (b) *and* $B_k$ *are independent* Bernoulli($\alpha_k$) *random variables for* $k = k_0 + 1, \ldots, K$.

REMARK 9. Although the result (b) requires an assumption of orthogonality, this restriction can be removed. See equation (38) of Corollary 4.5 from Leeb and Pötscher (2003).

Theorem 8 shows that forward stepwise is the best procedure in orthogonal designs. Suppose that $\alpha_k = \alpha > 0$ for each $k$. Then the limiting probability of correctly estimating $k_0$ is $\mathbb{P}\{k_F = k_0\} = (1 - \alpha)$ for forward stepwise, while for OLS-hard and backward stepwise, it is $(1-\alpha)^{K-k_0}$. Notice if $K - k_0$ is large, this last probability is approximated by $\exp(-(K - k_0)\alpha)$, which becomes exponentially small as $K$ increases. Simply put, the OLS-hard and backward stepwise methods are prone to overfitting. Figure 8 illustrates how the limiting probabilities vary under various choices for $K$ and $k_0$ (all figures computed with $\alpha = 0.10$). One can clearly see the superiority of the forward procedure, especially as $K$ increases.

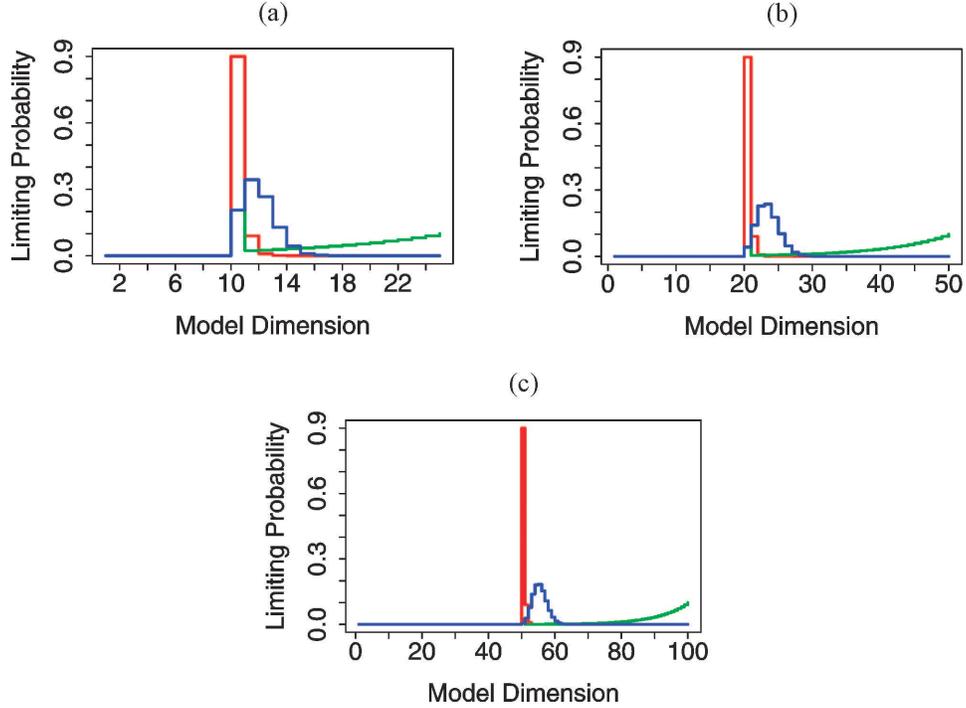

Fig. 8. *Complexity recovery in the orthogonal case. Limiting probabilities versus model dimension $k$ for the three estimators $\hat{k}_F$ (—), $\hat{k}_B$ (—) and $\hat{k}_O$ (—): (a) $K = 25$, $k_0 = 10$, (b) $K = 50$, $k_0 = 20$, (c) $K = 100$, $k_0 = 50$. In all cases $\alpha_k = 0.10$.*

**7. Diabetes data example.** As an illustration of the different model selection procedures we consider an example from Efron, Hastie, Johnstone and Tibshirani (2004). In illustrating the LARS method, Efron, Hastie, Johnstone and Tibshirani analyzed a diabetes study involving $n = 442$ patients in which the response of interest, $Y_i$, is a quantitative measure of disease progression recorded one year after baseline measurement. Data included ten baseline variables: age, sex, body mass index, average blood pressure and six blood serum measurements. All covariates were standardized and $Y_i$ was centered so that its mean value was zero. Two linear regression models were considered in the paper. The first was a main effects model involving the 10 baseline measurements, the second, a "quadratic model," which we re-analyze here, was made up of 64 covariates containing the 10 baseline measurements, 45 interactions for the 10 original covariates and 9 squared terms (these being the squares of each of the original covariates except for the gender variable, which is binary).

Table 1 contains the results from our analysis of the quadratic model. Listed are the top 10 variables as ranked by their absolute posterior means, $|\widehat{\beta}^*_{k,n}|$. Using an $\alpha = 0.10$ criteria, Zcut chooses a model with six variables starting from the top variable "bmi" (body mass index) and ending with "age.sex" (the age–sex interaction effect). The seventh variable, "bmi.map" (the interaction of body mass index and map, a blood pressure measurement), is borderline significant. Table 1 also reports results using OLS-hard, svsForwd and a new procedure, "OLSForwd" (all using an $\alpha = 0.10$ value). OLSForwd is the direct analogue of svsForwd, but orders $\boldsymbol{\beta}$ using $Z$-statistics $\widehat{Z}_{k,n}$ in place of the posterior mean. For all procedures the values in Table 1 are $Z$-statistics (12) derived from the restricted OLS for the selected model. This was done to allow direct comparison to the posterior mean values recorded in column 2.

Table 1 shows that the OLS-hard model differs significantly from Zcut. It excludes both "ltg" and "hdl" (blood serum measurements), both of which have large posterior mean values. We are not confident in the OLS-hard and suspect it is missing true signal here. The same comment applies to OLSForwd, which has produced the same model as OLS-hard. Note how svsForwd, the counterpart for OLSForwd, agrees closely with Zcut (it disagrees only on bmi.map, which is borderline significant). We believe the SVS models are more accurate than the OLS ones. In the next section we more systematically study the differences between the four procedures.

REMARK 10. Figure 9 displays the posterior density for $\sigma^2$. Note how the posterior is concentrated near one. This is typical of what we see in practice.

TABLE 1
*Top* 10 *variables from diabetes data (ranking based on absolute posterior means $|\widehat{\beta}^*_{k,n}|$). Entries for model selection procedures are $Z$-statistics* (12) *derived from the restricted OLS for the selected model*

|    | Variable | $\widehat{\beta}^*_{k,n}$ | Zcut  | OLS-hard | svsForwd | OLSForwd |
|----|----------|-------------------------|-------|----------|----------|----------|
| 1  | bmi      | 9.54                    | 8.29  | 13.70    | 8.15     | 13.70    |
| 2  | ltg      | 9.25                    | 7.68  | 0.00     | 7.82     | 0.00     |
| 3  | map      | 5.64                    | 5.39  | 7.06     | 4.99     | 7.06     |
| 4  | hdl      | $-4.37$                 | $-4.20$ | 0.00   | $-4.31$  | 0.00     |
| 5  | sex      | $-3.38$                 | $-4.03$ | $-1.95$ | $-4.02$ | $-1.95$  |
| 6  | age.sex  | 2.43                    | 3.58  | 3.19     | 3.47     | 3.19     |
| 7  | bmi.map  | 1.61                    | 0.00  | 2.56     | 3.28     | 2.56     |
| 8  | glu.2    | 0.84                    | 0.00  | 0.00     | 0.00     | 0.00     |
| 9  | bmi.2    | 0.46                    | 0.00  | 0.00     | 0.00     | 0.00     |
| 10 | tc.tch   | $-0.44$                 | 0.00  | 0.00     | 0.00     | 0.00     |

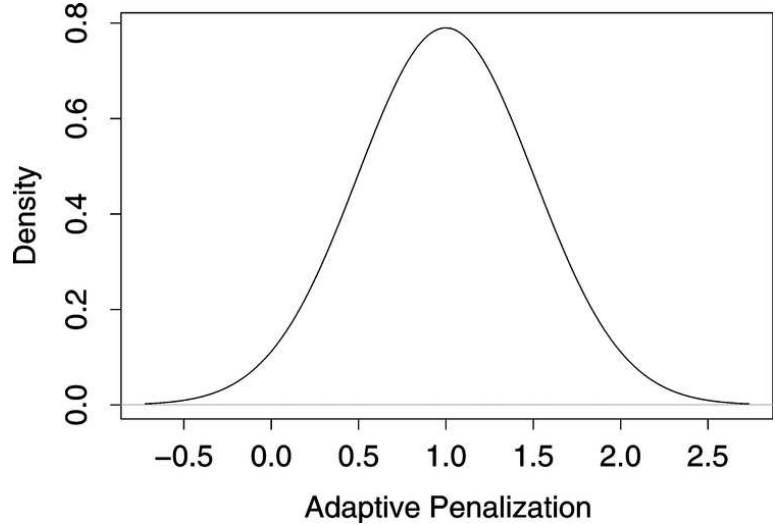

Fig. 9. *Posterior density for $\sigma^2$ from diabetes analysis.*

**8. Breiman simulations.** We used simulations to more systematically study performance. These followed the recipe given in Breiman (1992). Specifically, data were generated by taking $\varepsilon_i$ to be i.i.d. $N(0,1)$ variables, while covariates $\mathbf{x}_i$ were simulated independently from a multivariate normal distribution such that $\mathbb{E}(\mathbf{x}_i) = \mathbf{0}$ and $\mathbb{E}(x_{i,j} x_{i,k}) = \rho^{|j-k|}$, where $0 < \rho < 1$ was a correlation parameter. We considered two settings for $\rho$: (i) an uncorrelated design, $\rho = 0$; (ii) a correlated design, $\rho = 0.90$. For each $\rho$ setting we also considered two different sample size and model dimension configurations: (A) $n = 200$ and $K = 100$; (B) $n = 800$ and $K = 400$. Note that our illustrative example of Figure 1 corresponds to the Monte Carlo experiment (B) with $\rho = 0$.

In the higher-dimensional simulations (B), the nonzero $\beta_{k,0}$ coefficients were in 15 clusters of 7 adjacent variables centered at every 25th variable. For example, for the variables clustered around the 25th variable, the coefficient values were given by $\beta_{25+j,0} = |h - j|^{1.25}$ for $|j| < h$, where $h = 4$. The other 14 clusters were defined similarly. All other coefficients were set to zero. This gave a total of 105 nonzero values and 295 zero values. Coefficient values were adjusted by multiplying by a common constant to make the theoretical $R^2$ value equal to 0.75 [see Breiman (1992) for a discussion of this point].

Simulations (B) reflect a regression framework with a large number of zero coefficients. In contrast, simulations (A) were designed to represent a regression model with many weakly informative covariates. For (A), nonzero $\beta_{k,0}$ coefficients were grouped into 9 clusters each of size 5 centered at every 10th variable. Each of the 45 nonzero coefficients was set to the same value. Coefficient values were then adjusted by multiplying by a common constant to make the theoretical $R^2$ value equal to 0.75. This ensured that the overall signal to noise ratio was the same as (B), but with each coefficient having less explanatory power.

Simulations were repeated 100 times independently for each of the four experiments. Results are recorded in Table 2 for each of the procedures Zcut, svsForwd, OLS-hard and OLSForwd (all using an $\alpha = 0.10$ value). Table 2 records what we call "TotalMiss," "FDR" and "FNR." The TotalMiss is the total number of misclassified variables, that is, the total number of falsely identified nonzero $\beta_{k,0}$ coefficients and falsely identified zero coefficients. This is an unbiased estimator for the risk discussed in Theorem 5. The FDR and FNR are the false discovery and false nondiscovery rates defined as the false positive and false negative rates for those coefficients identified as nonzero and zero, respectively. The TotalMiss, FDR and FNR values reported are the averaged values from the 100 simulations. Also recorded is $\hat{k}$, the average number of variables selected by a procedure. Table 2 also includes the performance value "Perf," a measure of prediction accuracy, defined as

$$\text{Perf} = 1 - \frac{\|\mathbf{X}\widehat{\boldsymbol{\beta}} - \mathbf{X}\boldsymbol{\beta}_0\|^2}{\|\mathbf{X}\boldsymbol{\beta}_0\|^2},$$

where $\widehat{\boldsymbol{\beta}}$ is the estimator for $\boldsymbol{\beta}_0$. So Perf equals zero when $\widehat{\boldsymbol{\beta}} = \mathbf{0}$ and equals one when $\widehat{\boldsymbol{\beta}} = \boldsymbol{\beta}_0$. The value for Perf was again averaged over the 100 simulations.

TABLE 2
*Breiman simulations*

|  | $\rho = 0$ (uncorrelated X) | | | | | $\rho = 0.9$ (correlated X) | | | | |
|---|---|---|---|---|---|---|---|---|---|---|
|  | $\hat{k}$ | Perf | TotalMiss | FDR | FNR | $\hat{k}$ | Perf | TotalMiss | FDR | FNR |
| (A) *Moderate number of covariates with few* (55%) *that are zero* ($n = 200$, $K = 100$ *and 55 zero* $\beta_{k,0}$). | | | | | | | | | | |
| Zcut | 41.44 | 0.815 | 11.99 | 0.097 | 0.129 | 10.06 | 0.853 | 38.49 | 0.167 | 0.408 |
| svsForwd | 34.02 | 0.753 | 15.09 | 0.054 | 0.191 | 8.31 | 0.826 | 39.39 | 0.156 | 0.415 |
| OLS-hard | 41.99 | 0.791 | 14.06 | 0.128 | 0.145 | 11.08 | 0.707 | 45.31 | 0.496 | 0.446 |
| OLSForwd | 26.90 | 0.612 | 20.92 | 0.042 | 0.258 | 5.96 | 0.574 | 44.64 | 0.459 | 0.445 |
| (B) *Large number of covariates with many* (74%) *that are zero* ($n = 800$, $K = 400$ *and 295 zero* $\beta_{k,0}$). | | | | | | | | | | |
| Zcut | 75.96 | 0.903 | 39.62 | 0.068 | 0.106 | 36.67 | 0.953 | 72.61 | 0.055 | 0.194 |
| svsForwd | 86.81 | 0.904 | 41.19 | 0.130 | 0.095 | 24.42 | 0.926 | 81.90 | 0.025 | 0.216 |
| OLS-hard | 106.74 | 0.883 | 58.54 | 0.279 | 0.097 | 45.41 | 0.706 | 121.37 | 0.676 | 0.255 |
| OLSForwd | 61.09 | 0.846 | 49.87 | 0.046 | 0.138 | 9.14 | 0.303 | 106.48 | 0.590 | 0.259 |

REMARK 11. Given the high dimensionality of the simulations, both svsForwd and OLSForwd often stopped early and produced models that were much too small. To compensate, we slightly altered their definitions. For svsForwd, we modified the definition of $\hat{k}_F$ [cf. (16)] to

$$\hat{k}_F = \min\{k-1 : |\widetilde{Z}_{k,n}| < z_{\alpha_k/2} \text{ and } |\widehat{\beta}^*_{k,n}| \leq C, k = 1, \ldots, K+1\},$$

where $C = 3$. In this way, svsForwd stops the first time the null hypothesis is not rejected *and* if the absolute posterior mean is no longer a large value. The definition for OLSForwd was altered in similar fashion, but using $\widehat{Z}_{k,n}$ in place of $\widehat{\beta}^*_{k,n}$.

8.1. *Results.* The simulations revealed several interesting patterns, summarized as follows:

1. Zcut beats OLS-hard across all performance categories. It maintains low risk, has small FDR values and has good prediction error performance in both the near-orthogonal (uncorrelated) and nonorthogonal (correlated) **X** cases. Performance differences between Zcut and OLS-hard become more appreciable in the near-orthogonal simulation (B) involving many zero coefficients, because this is when the effect of selective shrinkage is most pronounced. For example, the OLS-hard misclassifies about 19 coefficients more on average, and has a FDR more than 4 times larger than Zcut's. Large gains are also seen in the correlated case (B). There, the OLS-hard misclassifies over 48 more coefficients on average than Zcut and its FDR is more than 12 times higher.
2. It is immediately clear upon comparing svsForwd to OLSForwd that SVS is capable of some serious model averaging. These two procedures differ only in the way they rank coefficients, so the disparity in their two performances is clear evidence of SVS's ability to model average.
3. In the $\rho = 0$ simulations, svsForwd is roughly the same as OLS-hard in simulation (A) and significantly better in simulation (B). In the correlated setting, svsForwd is significantly better. Thus, overall svsForwd is as good, and in most cases significantly better, than OLS-hard. This suggests that svsForwd is starting to tap into the oracle property forward stepwise has relative to OLS-hard and provides indirect evidence that SVS is capable of reducing model uncertainty in finite samples.
4. It is interesting to note how badly OLSForwd performs relative to OLS-hard in simulation (A) when $\rho = 0$. In orthogonal designs, OLSForwd is equivalent to OLS-hard, but the $\rho = 0$ design is only near-orthogonal. With only a slight departure from orthogonality, we see the importance of a reliable ranking for the coordinates of $\boldsymbol{\beta}$. Note that this effect is less pronounced in simulation (B) because of the larger sample size. This is because $\mathbf{X}^t\mathbf{X}/n \overset{\text{a.s.}}{\to} \mathbf{I}$ as $n \to \infty$, so simulation (B) should be closer to orthogonality.

5. While our theory does not cover Zcut's performance in correlated settings, it is interesting to note how well it does in the $\rho = 0.9$ simulations relative to OLS-hard. The explanation for its success here, however, is probably different from that for the orthogonal setting. For example, it is possible that its performance gains may be mostly due to the use of generalized ridge estimators. As is well known, such estimators are much more stable than OLS in multicollinear settings. We should also note that while Zcut is better than OLS-hard here, its performance relative to the orthogonal simulations is noticeably worse. This is not unexpected though. Correlation has the effect of reducing the dimension of the problem. So performance measurements like TotalMiss and FDR will naturally be less favorable.

## APPENDIX: PROOFS

PROOF OF THEOREM 2. We start by establishing that $\widehat{\boldsymbol{\beta}}_n^\circ$ is consistent, which is part of the conclusion of Theorem 2. First observe that

$$\widehat{\boldsymbol{\beta}}_n^\circ = n^{-1}\boldsymbol{\Sigma}_n^{-1}\mathbf{X}^t\mathbf{Y} = \boldsymbol{\beta}_0 + \boldsymbol{\Delta}_n,$$

where $\boldsymbol{\Delta}_n = \boldsymbol{\Sigma}_n^{-1}\mathbf{X}^t\boldsymbol{\epsilon}/n$ and $\boldsymbol{\epsilon} = (\varepsilon_1, \ldots, \varepsilon_n)^t$. From $\mathbb{E}(\boldsymbol{\Delta}_n) = \mathbf{0}$ and $\mathrm{Var}(\boldsymbol{\Delta}_n) = \sigma_0^2 \boldsymbol{\Sigma}_n^{-1}/n$, it is clear that $\widehat{\boldsymbol{\beta}}_n^\circ \xrightarrow{\mathrm{p}} \boldsymbol{\beta}_0$. Next, a little bit of rearrangement shows that

$$\widehat{\boldsymbol{\theta}}_n^*(\boldsymbol{\gamma}, \sigma^2) = (\mathbf{I} - (\sigma^{-2}\lambda_n^{-1}\mathbf{X}^t\mathbf{X} + \boldsymbol{\Gamma}^{-1})^{-1}\boldsymbol{\Gamma}^{-1})\widehat{\boldsymbol{\beta}}_n^\circ.$$

Consequently,

$$\widehat{\boldsymbol{\theta}}_n^* = \widehat{\boldsymbol{\beta}}_n^\circ - \int (\sigma^{-2}\lambda_n^{-1}\mathbf{X}^t\mathbf{X} + \boldsymbol{\Gamma}^{-1})^{-1}\boldsymbol{\Gamma}^{-1}\widehat{\boldsymbol{\beta}}_n^\circ (\pi \times \mu)(d\boldsymbol{\gamma}, d\sigma^2 | \mathbf{Y}^*)$$

$$= \widehat{\boldsymbol{\beta}}_n^\circ - \lambda_n^* \int \sigma^2 \mathbf{V}_n^{-1}\boldsymbol{\Gamma}^{-1}\widehat{\boldsymbol{\beta}}_n^\circ (\pi \times \mu)(d\boldsymbol{\gamma}, d\sigma^2 | \mathbf{Y}^*),$$

where $\lambda_n^* = \lambda_n/n$ and $\mathbf{V}_n = \boldsymbol{\Sigma}_n + \sigma^2\lambda_n^*\boldsymbol{\Gamma}^{-1}$. By the Jordan decomposition theorem, we can write $\mathbf{V}_n = \sum_{k=1}^K e_{k,n}\mathbf{v}_{k,n}\mathbf{v}_{k,n}^t$, where $\{\mathbf{v}_{k,n}\}$ is a set of orthonormal eigenvectors with eigenvalues $\{e_{k,n}\}$. For convenience, assume that the eigenvalues have been ordered so that $e_{1,n} \leq \cdots \leq e_{K,n}$. The assumption that $\boldsymbol{\Sigma}_n \to \boldsymbol{\Sigma}_0$, where $\boldsymbol{\Sigma}_0$ is positive definite, ensures that the minimum eigenvalue for $\boldsymbol{\Sigma}_n$ is larger than some $e_0 > 0$ for sufficiently large $n$. Therefore, if $n$ is large enough,

$$e_{1,n} \geq e_0 + \sigma^2\lambda_n^* \min_k \gamma_k^{-1} \geq e_0 > 0.$$

Notice that

$$\|\mathbf{V}_n^{-1}\boldsymbol{\Gamma}^{-1}\widehat{\boldsymbol{\beta}}_n^\circ\|^2 = \sum_{k=1}^K e_{k,n}^{-2}(\mathbf{v}_{k,n}^t\boldsymbol{\Gamma}^{-1}\widehat{\boldsymbol{\beta}}_n^\circ)^2 \leq e_0^{-2}\|\widehat{\boldsymbol{\beta}}_n^\circ\|^2 \sum_{k=1}^K \gamma_k^{-2}.$$

Thus, since $\gamma_k \geq \eta_0$ over the support of $\pi$, and $\sigma^2 \leq s_0^2$ over the support of $\mu$,

$$\left\| \int \sigma^2 \mathbf{V}_n^{-1} \boldsymbol{\Gamma}^{-1} \widehat{\boldsymbol{\beta}}_n^\circ (\pi \times \mu)(d\boldsymbol{\gamma}, d\sigma^2 | \mathbf{Y}^*) \right\|$$

$$\leq e_0^{-1} \|\widehat{\boldsymbol{\beta}}_n^\circ\| \left( \sum_{k=1}^K \int \sigma^4 \gamma_k^{-2} (\pi \times \mu)(d\boldsymbol{\gamma}, d\sigma^2 | \mathbf{Y}^*) \right)^{1/2}$$

$$\leq \frac{K^{1/2} s_0^2}{\eta_0 e_0} \|\widehat{\boldsymbol{\beta}}_n^\circ\|.$$

Deduce that $\widehat{\boldsymbol{\theta}}_n^* = \widehat{\boldsymbol{\beta}}_n^\circ + O_p(\lambda_n^*) \xrightarrow{\text{p}} \boldsymbol{\beta}_0$. $\square$

Before proving Theorem 3, we state a lemma. This will also be useful in the proofs of some later theorems.

LEMMA A.1. *Assume that for each $n$, $\varepsilon_{n1}, \ldots, \varepsilon_{nn}$ are independent random variables such that $\mathbb{E}(\varepsilon_{ni}) = 0$, $\mathbb{E}(\varepsilon_{ni}^2) = \sigma_0^2$ and $\mathbb{E}(\varepsilon_{ni}^4) \leq M$ for some finite $M$. If* (D1)–(D4) *hold, then*

$$(18) \qquad n^{-1/2} \mathbf{X}^t \boldsymbol{\epsilon}_n = n^{-1/2} \sum_{i=1}^n \varepsilon_{ni} \mathbf{x}_i \xrightarrow{\text{d}} \mathrm{N}(\mathbf{0}, \sigma_0^2 \boldsymbol{\Sigma}_0),$$

*where $\boldsymbol{\epsilon}_n = (\varepsilon_{n1}, \ldots, \varepsilon_{nn})^t$.*

PROOF. Let $S_n = \sum_{i=1}^n \varepsilon_{ni} \mathbf{x}_i^t \boldsymbol{\ell} / \sqrt{n}$, where $\boldsymbol{\ell} \in \mathbb{R}^K$ is some arbitrary nonzero vector. Let $s_n^2 = \sigma_0^2 \boldsymbol{\ell}^t \boldsymbol{\Sigma}_n \boldsymbol{\ell}$ and define $\zeta_{ni} = n^{-1/2} \varepsilon_{ni} \mathbf{x}_i^t \boldsymbol{\ell} / s_n$. Then, $S_n/s_n = \sum_{i=1}^n \zeta_{ni}$, where $\zeta_{ni}$ are independent random variables such that $\mathbb{E}(\zeta_{ni}) = 0$ and $\sum_{i=1}^n \mathbb{E}(\zeta_{ni}^2) = 1$. To prove (18), we will verify the Lindeberg condition

$$\sum_{i=1}^n \mathbb{E}(\zeta_{ni}^2 \mathbb{I}\{|\zeta_{ni}| \geq \delta\}) \to 0 \qquad \text{for each } \delta > 0,$$

where $\mathbb{I}(\cdot)$ denotes the indicator function. This will show that $S_n/s_n \xrightarrow{\text{d}} \mathrm{N}(0,1)$, which by the Cramér–Wold device implies (18) because $s_n^2 \to \sigma_0^2 \boldsymbol{\ell}^t \boldsymbol{\Sigma}_0 \boldsymbol{\ell}$. Observe that

$$\mathbb{E}(\zeta_{ni}^2 \mathbb{I}\{|\zeta_{ni}| > \delta\}) = \frac{(\mathbf{x}_i^t \boldsymbol{\ell})^2}{n s_n^2} \mathbb{E}(\varepsilon_{ni}^2 \mathbb{I}\{|\varepsilon_{ni}| \geq r_{ni} s_n \delta\}),$$

where $r_{ni} = \sqrt{n}/|\mathbf{x}_i^t \boldsymbol{\ell}|$. By the Cauchy–Schwarz inequality and the assumption of a bounded fourth moment for $\varepsilon_{ni}$,

$$\mathbb{E}(\varepsilon_{ni}^2 \mathbb{I}\{|\varepsilon_{ni}| \geq r_{ni} s_n \delta\}) \leq (\mathbb{E}(\varepsilon_{ni}^4) \mathbb{P}\{|\varepsilon_{ni}| \geq r_{ni} s_n \delta\})^{1/2} \leq \frac{M^{1/2} \sigma_0}{r_{ni} s_n \delta}.$$

Bound $r_{ni}$ below by

$$r_{ni} \geq r_n := \left(\max_{1 \leq i \leq n} |\mathbf{x}_i^t \boldsymbol{\ell}|/\sqrt{n}\right)^{-1}.$$

Notice that $r_n \to \infty$ by the assumption that $\max_i \|\mathbf{x}_i\|/\sqrt{n} \to 0$. Substituting the bound for $r_{ni}$, and since $(\mathbf{x}_i^t \boldsymbol{\ell})^2$ sums to $ns_n^2/\sigma_0^2$, and $s_n^2$ remains bounded away from zero since $s_n^2 \to \sigma_0^2 \boldsymbol{\ell}^t \boldsymbol{\Sigma}_0 \boldsymbol{\ell}$,

$$\sum_{i=1}^n \mathbb{E}(\zeta_{ni}^2 \mathbb{I}\{|\zeta_{ni}| > \delta\}) \leq \frac{M^{1/2}}{\sigma_0 r_n s_n \delta} \to 0. \qquad \square$$

PROOF OF THEOREM 3. Let $\phi(\cdot|m, \tau^2)$ denote a normal density with mean $m$ and variance $\tau^2$. By dividing the numerator and denominator by $n^{-K/2} \prod_{i=1}^n \phi(Y_{ni}^*|\mathbf{x}_i^t \boldsymbol{\beta}_0, n)$, one can show that

$$(19) \quad \frac{\nu_n(S(\boldsymbol{\beta}_1, C/\sqrt{n})|\mathbf{Y}_n^*)}{\nu_n(S(\boldsymbol{\beta}_0, C/\sqrt{n})|\mathbf{Y}_n^*)} = \frac{n^{K/2} \int \mathbb{I}\{\boldsymbol{\beta} \in S(\boldsymbol{\beta}_1, C/\sqrt{n})\} L_n(\boldsymbol{\beta}) \nu(d\boldsymbol{\beta})}{n^{K/2} \int \mathbb{I}\{\boldsymbol{\beta} \in S(\boldsymbol{\beta}_0, C/\sqrt{n})\} L_n(\boldsymbol{\beta}) \nu(d\boldsymbol{\beta})},$$

where

$$\log(L_n(\boldsymbol{\beta})) = -\tfrac{1}{2}(\boldsymbol{\beta} - \boldsymbol{\beta}_0)^t \boldsymbol{\Sigma}_n (\boldsymbol{\beta} - \boldsymbol{\beta}_0) + n^{-1/2} \sum_{i=1}^n \varepsilon_{ni} \mathbf{x}_i^t (\boldsymbol{\beta} - \boldsymbol{\beta}_0).$$

Consider the denominator in (19). By making the change of variables from $\boldsymbol{\beta}$ to $\mathbf{u} = \sqrt{n}(\boldsymbol{\beta} - \boldsymbol{\beta}_0)$, we can rewrite this as

$$(20) \quad \int \mathbb{I}\{\mathbf{u} \in S(\mathbf{0}, C)\} L_{n0}(\mathbf{u}) f(\boldsymbol{\beta}_0 + n^{-1/2} \mathbf{u}) \, d\mathbf{u},$$

where $\log(L_{n0}(\mathbf{u})) = g_n(\mathbf{u}) + O(1/n)$ and $g_n(\mathbf{u}) = \sum_{i=1}^n \varepsilon_{ni} \mathbf{x}_i^t \mathbf{u}/n$. The $O(1/n)$ term corresponds to $\mathbf{u}^t \boldsymbol{\Sigma}_n \mathbf{u}/n$ and is uniform over $\mathbf{u} \in S(\mathbf{0}, C)$. Observe that, for each $\delta > 0$,

$$\mathbb{P}\{|g_n(\mathbf{u})| \geq \delta\} \leq \frac{1}{\delta^2 n^2} \sum_{i=1}^n \mathbb{E}(\varepsilon_{ni} \mathbf{x}_i^t \mathbf{u})^2 \leq \frac{C^2}{\delta^2 n^2} \sum_{i=1}^n \|\mathbf{x}_i\|^2 = o(1),$$

where the last inequality on the right-hand side follows from the Cauchy–Schwarz inequality and from the assumption that $\max_i \|\mathbf{x}_i\|/\sqrt{n} = o(1)$. Therefore, $g_n(\mathbf{u}) \xrightarrow{P} 0$ uniformly over $\mathbf{u} \in S(\mathbf{0}, C)$. Because $f$ is continuous [and keeping in mind it remains positive and bounded over $S(\mathbf{0}, C)$], deduce that the log of (20) converges in probability to

$$(21) \quad \log(f(\boldsymbol{\beta}_0)) + \log\left(\int \mathbb{I}\{\mathbf{u} \in S(\mathbf{0}, C)\} \, d\mathbf{u}\right).$$

Meanwhile, for the numerator in (19), make the change of variables from $\boldsymbol{\beta}$ to $\mathbf{u} = \sqrt{n}(\boldsymbol{\beta} - \boldsymbol{\beta}_1)$ to rewrite this as

$$(22) \quad \int \mathbb{I}\{\mathbf{u} \in S(\mathbf{0}, C)\} L_{n1}(\mathbf{u}) f(\boldsymbol{\beta}_1 + n^{-1/2} \mathbf{u}) \, d\mathbf{u},$$

where

$$\log(L_{n1}(\mathbf{u})) = -\tfrac{1}{2}(\boldsymbol{\beta}_1 - \boldsymbol{\beta}_0)^t \boldsymbol{\Sigma}_0 (\boldsymbol{\beta}_1 - \boldsymbol{\beta}_0)$$
$$+ n^{-1/2} \sum_{i=1}^{n} \varepsilon_{ni} \mathbf{x}_i^t (\boldsymbol{\beta}_1 - \boldsymbol{\beta}_0) + g_n(\mathbf{u}) + o(1)$$

uniformly over $\mathbf{u} \in S(\mathbf{0}, C)$. Consider the second term on the right-hand side of the last expression. Set $\boldsymbol{\ell} = \boldsymbol{\beta}_1 - \boldsymbol{\beta}_0$. By Lemma A.1, since $\sigma_0^2 = 1$, it follows that

$$n^{-1/2} \sum_{i=1}^{n} \varepsilon_{ni} \mathbf{x}_i^t \boldsymbol{\ell} \xrightarrow{d} \mathrm{N}(0, \boldsymbol{\ell}^t \boldsymbol{\Sigma}_0 \boldsymbol{\ell}).$$

Now extract the expressions not depending upon $\mathbf{u}$ outside the integral in (22), take logs and use $g_n(\mathbf{u}) \xrightarrow{P} 0$ uniformly over $\mathbf{u} \in S(\mathbf{0}, C)$ to deduce that the log of (22) converges in distribution to

$$\begin{aligned}(23)\quad & -\tfrac{1}{2}(\boldsymbol{\beta}_1 - \boldsymbol{\beta}_0)^t \boldsymbol{\Sigma}_0 (\boldsymbol{\beta}_1 - \boldsymbol{\beta}_0) \\ & + (\boldsymbol{\beta}_1 - \boldsymbol{\beta}_0)^t \mathbf{Z} + \log(f(\boldsymbol{\beta}_1)) + \log\left(\int \mathbb{I}\{\mathbf{u} \in S(\mathbf{0}, C)\}\, d\mathbf{u}\right).\end{aligned}$$

To complete the proof, take the difference of (23) and (21) and note the cancellation of the logs of $\int \mathbb{I}\{\mathbf{u} \in S(\mathbf{0}, C)\}\, d\mathbf{u}$. $\square$

PROOF OF THEOREM 4. First note that

$$\widehat{\boldsymbol{\beta}}_{nn}^*(\boldsymbol{\gamma}_0) = (\mathbf{X}^t \mathbf{X} + n \boldsymbol{\Gamma}_0^{-1})^{-1} \mathbf{X}^t \mathbf{Y}_n^* = (\boldsymbol{\Sigma}_0 + \boldsymbol{\Gamma}_0^{-1})^{-1} \boldsymbol{\Sigma}_0 \boldsymbol{\beta}_0 + n^{-1/2} \mathbf{V}_n^{-1} \mathbf{X}^t \boldsymbol{\epsilon}_n + o(1),$$

where $\mathbf{V}_n = \boldsymbol{\Sigma}_n + \boldsymbol{\Gamma}_0^{-1}$. The $o(1)$ term on the right-hand side is due to $\boldsymbol{\Sigma}_n \to \boldsymbol{\Sigma}_0$. Also, by Lemma A.1, the second term on the right-hand side converges in distribution to $(\boldsymbol{\Sigma}_0 + \boldsymbol{\Gamma}_0^{-1})^{-1} \mathbf{Z}$, where $\mathbf{Z}$ has a $\mathrm{N}(\mathbf{0}, \boldsymbol{\Sigma}_0)$ distribution. Deduce that $\widehat{\boldsymbol{\beta}}_{nn}^*(\boldsymbol{\gamma}_0) \xrightarrow{d} Q(\cdot | \boldsymbol{\gamma}_0)$. $\square$

PROOF OF THEOREM 5. A little algebra (keeping in mind $\boldsymbol{\Sigma}_n = \mathbf{I}$) shows $\widehat{\boldsymbol{\beta}}_n^* = \sqrt{n}(\mathbf{I} + \boldsymbol{\Gamma}_0^{-1})^{-1} \widehat{\boldsymbol{\beta}}_n^\circ / \widehat{\sigma}_n$. Hence, recalling the definition (12) for $\widehat{Z}_{k,n}$,

$$\widehat{\beta}_{k,n}^* = d_{k,0} \times \frac{n^{1/2} \widehat{\beta}_{k,n}^\circ}{\widehat{\sigma}_n} = d_{k,0} \widehat{Z}_{k,n},$$

where $d_{k,0} = \gamma_{k,0}/(1 + \gamma_{k,0})$ and the last equality holds because $s_{kk} = 1$. Under the assumption of normality, $\sqrt{n}\widehat{\beta}_{k,n}^\circ$ has a $\mathrm{N}(m_{k,n}, \sigma_0^2)$ distribution, where $m_{k,n} = \sqrt{n}\beta_{k,0}$. Choose $\boldsymbol{\gamma}_0$ such that $d_{k,0} = \delta_1$ for each $k \in \mathscr{B}_0$ and $d_{k,0} = \delta_2$ for each $k \in \mathscr{B}_0^c$, where $0 < \delta_1, \delta_2 < 1$ are values to be specified.

Therefore,

$$\begin{aligned}
&\mathscr{R}_O(\alpha) - \mathscr{R}_Z(\alpha) \\
&= (K - k_0)(\mathbb{P}\{|N(0, \sigma_0^2)| \geq \widehat{\sigma}_n z_{\alpha/2}\} - \mathbb{P}\{|N(0, \sigma_0^2)| \geq \delta_1^{-1}\widehat{\sigma}_n z_{\alpha/2}\}) \\
&\quad + \sum_{k \in \mathscr{B}_0^c} (\mathbb{P}\{|N(m_{k,n}, \sigma_0^2)| < \widehat{\sigma}_n z_{\alpha/2}\} - \mathbb{P}\{|N(m_{k,n}, \sigma_0^2)| < \delta_2^{-1}\widehat{\sigma}_n z_{\alpha/2}\}),
\end{aligned}$$

where the $\mathbb{P}$-distributions on the right-hand side correspond to the joint distribution for a normal random variable and the distribution for $\widehat{\sigma}_n$, where $\widehat{\sigma}_n^2/\sigma_0^2$ has an independent $\chi^2$-distribution with $n - K$ degrees of freedom. It is clear that the sum on the right-hand side can be made arbitrarily close to zero, uniformly for $\alpha \in [\delta, 1 - \delta]$, by choosing $\delta_2$ close to one, while the first term on the right-hand side remains positive and uniformly bounded away from zero over $\alpha \in [\delta, 1 - \delta]$ whatever the choice for $\delta_1$. Thus, for a suitably chosen $\delta_2$, $\mathscr{R}_O(\alpha) - \mathscr{R}_Z(\alpha) > 0$ for each $\alpha \in [\delta, 1 - \delta]$. □

PROOF OF THEOREM 6. Choose some $j \in \mathscr{B}_0^c$. Let $A_j = \{\boldsymbol{\gamma}: d_j \leq 1 - \delta\}$, where $d_j = \gamma_j/(1 + \gamma_j)$. To prove part (a), we show that

$$\pi_n(A_j|\mathbf{Y}^*) = \frac{\int_{A_j} f(\mathbf{Y}^*|\boldsymbol{\gamma})\pi(d\boldsymbol{\gamma})}{\int f(\mathbf{Y}^*|\boldsymbol{\gamma})\pi(d\boldsymbol{\gamma})} \xrightarrow{\mathrm{P}} 0.$$

By definition, $f(\mathbf{Y}^*|\boldsymbol{\gamma}) = \int f(\mathbf{Y}^*|\boldsymbol{\beta})f(\boldsymbol{\beta}|\boldsymbol{\gamma})\,d\boldsymbol{\beta}$, where

$$f(\mathbf{Y}^*|\boldsymbol{\beta})f(\boldsymbol{\beta}|\boldsymbol{\gamma}) = C\exp\Big(-\frac{1}{2n}(\mathbf{Y}^* - \mathbf{X}\boldsymbol{\beta})^t(\mathbf{Y}^* - \mathbf{X}\boldsymbol{\beta}) - \frac{1}{2}\boldsymbol{\beta}^t\boldsymbol{\Gamma}^{-1}\boldsymbol{\beta}\Big)|\boldsymbol{\Gamma}|^{-1/2}$$

and $C$ is a generic constant not depending upon $\boldsymbol{\gamma}$. By some straightforward calculations that exploit conjugacy and orthogonality,

$$(24) \qquad f(\mathbf{Y}^*|\boldsymbol{\gamma}) = C\exp\Big(\tfrac{1}{2}\sum_{k=1}^K d_k \xi_{k,n}^2\Big)\prod_{k=1}^K (1 + \gamma_k)^{-1/2},$$

where $(\xi_{1,n}, \ldots, \xi_{K,n})^t = \widehat{\sigma}_n^{-1} n^{-1/2} \mathbf{X}^t \mathbf{Y}$.

Let $B = \{\boldsymbol{\gamma}: 1 - \delta_k \leq d_k \leq 1 - \delta_k/2, k = 1, \ldots, K\}$, where $0 < \delta_k < 1$ are small values that will be specified. Observe that

$$\pi_n(A_j|\mathbf{Y}^*) \leq \frac{\int_{A_j} f(\mathbf{Y}^*|\boldsymbol{\gamma})\pi(d\boldsymbol{\gamma})}{\int_B f(\mathbf{Y}^*|\boldsymbol{\gamma})\pi(d\boldsymbol{\gamma})}.$$

Over the set $A_j$ we have the upper bound

$$f(\mathbf{Y}^*|\boldsymbol{\gamma}) \leq C\exp\Big\{\tfrac{1}{2}\Big(\sum_{k \in \mathscr{B}_0} \xi_{k,n}^2 + \sum_{k \in \mathscr{B}_0^c - \{j\}} \xi_{k,n}^2 + (1 - \delta)\xi_{j,n}^2\Big)\Big\}$$

because $0 < d_k < 1$, while over $B$ we have the lower bound

$$f(\mathbf{Y}^*|\boldsymbol{\gamma}) \geq C \exp\left\{\frac{1}{2}\left(\sum_{k \in \mathscr{B}_0^c - \{j\}} (1-\delta_k)\xi_{k,n}^2 + (1-\delta_j)\xi_{j,n}^2\right)\right\} \prod_{k=1}^{K}\left(\frac{2}{\delta_k}\right)^{-1/2}.$$

An application of Lemma A.1 (which also applies to nontriangular arrays) shows

$$(\xi_{1,n}, \ldots, \xi_{K,n})^t = \widehat{\sigma}_n^{-1}(n^{1/2}\boldsymbol{\beta}_0 + O_p(1)).$$

Therefore,

$$\pi_n(A_j|\mathbf{Y}^*) \leq \exp\left\{O_p(1) + \frac{n}{2\widehat{\sigma}_n^2}\left(\sum_{k \in \mathscr{B}_0^c - \{j\}} \delta_k \beta_{k,0}^2 \right.\right.$$
(25)
$$\left.\left. + (\delta_j - \delta)\beta_{j,0}^2 + O_p(1/\sqrt{n})\right)\right\}\frac{\pi(A_j)}{\pi(B)}.$$

Choose $\delta_j < \delta$. It is clear we can find a set of values $\{\delta_k : k \neq j\}$ chosen small enough so that

$$\sum_{k \in \mathscr{B}_0^c - \{j\}} \delta_k \beta_{k,0}^2 + (\delta_j - \delta)\beta_{j,0}^2 < 0.$$

This ensures that the expression in the exponent of (25) converges to $-\infty$ in probability. Note for this result we assume $\widehat{\sigma}_n^2$ has a nonzero limit (we give a proof shortly that $\widehat{\sigma}_n^2 \xrightarrow{P} \sigma_0^2$). Therefore, since $\pi(B)$ must be strictly positive for small enough $\delta_k > 0$ (by our assumptions regarding the support for $\pi$), conclude from (25) that $\pi_n(A_j|\mathbf{Y}^*) \xrightarrow{P} 0$.

To prove part (b), let $f_k(\gamma_k|w)$ denote the density for $\gamma_k$ given $w$. From (4), it is seen that $f_k(\gamma_k|w) = (1-w)g_0(\gamma_k) + wg_1(\gamma_k)$. Therefore,

$$f_k^*(\gamma_k|w) \propto f(\mathbf{Y}|\boldsymbol{\gamma})f_k(\gamma_k|w)$$
$$\propto \exp(\tfrac{1}{2}d_k\xi_{k,n}^2)(1+\gamma_k)^{-1/2}f_k(\gamma_k|w),$$

which is the expression (13). Furthermore, by Lemma A.1 deduce that $\xi_{k,n}$ converges to a standard normal if $\beta_{k,0} = 0$ (we are using $\widehat{\sigma}_n^2 \xrightarrow{P} \sigma_0^2$, which still needs to be proven).

To complete the proof, we now show $\widehat{\sigma}_n^2$ is consistent. For this proof we do not assume orthogonality, only that $\boldsymbol{\Sigma}_n$ is positive definite (this generality will be useful for later proofs). Observe that $\widehat{\sigma}_n^2 = (\boldsymbol{\epsilon}^t\boldsymbol{\epsilon} - \boldsymbol{\epsilon}^t\mathbf{H}\boldsymbol{\epsilon})/(n-K)$, where $\mathbf{H} = \mathbf{X}(\mathbf{X}^t\mathbf{X})^{-1}\mathbf{X}^t$. It follows from Chebyshev's inequality using the moment assumptions on $\varepsilon_i$ that $\boldsymbol{\epsilon}^t\boldsymbol{\epsilon}/(n-K) \xrightarrow{P} \sigma_0^2$, while from Markov's inequality, for each $\delta > 0$,

$$\mathbb{P}\{\boldsymbol{\epsilon}^t\mathbf{H}\boldsymbol{\epsilon} \geq (n-K)\delta\} \leq \frac{\mathbb{E}(\boldsymbol{\epsilon}^t\mathbf{H}\boldsymbol{\epsilon})}{(n-K)\delta} = \frac{\text{Trace}(\mathbf{H}\mathbb{E}(\boldsymbol{\epsilon}\boldsymbol{\epsilon}^t))}{(n-K)\delta} = \frac{K\sigma_0^2}{(n-K)\delta} \to 0.$$

Deduce that $\widehat{\sigma}_n^2 \xrightarrow{\text{P}} \sigma_0^2$. □

PROOF OF THEOREM 7. Under the assumption of orthogonality, and using the fact that $\lambda_n = n$, it follows that

$$\widehat{\boldsymbol{\beta}}_n^*(\boldsymbol{\gamma}, \sigma^2) = \widehat{\sigma}_n^{-1} n^{1/2} \mathbf{D} \boldsymbol{\beta}_0 + \widehat{\sigma}_n^{-1} n^{-1/2} \mathbf{D} \mathbf{X}^t \boldsymbol{\epsilon},$$

where $\mathbf{D}$ is the diagonal matrix $\text{diag}(d_1, \ldots, d_K)$ and $d_k = \gamma_k/(\gamma_k + \sigma^2)$. Taking expectations with respect to the posterior, deduce that

$$\widehat{\beta}_{k,n}^* = \widehat{\sigma}_n^{-1} n^{1/2} d_k^* \beta_{k,0} + \widehat{\sigma}_n^{-1} d_k^* \zeta_{k,n}, \tag{26}$$

where $d_k^* = \mathbb{E}(d_k | \mathbf{Y}^*)$ and $\zeta_{k,n}$ is the $k$th coordinate of $\mathbf{X}^t \boldsymbol{\epsilon}/\sqrt{n}$. From Lemma A.1 and $\widehat{\sigma}_n^2 \xrightarrow{\text{P}} \sigma_0^2$ (proven in Theorem 6), we have $\widehat{\sigma}_n^{-1} \mathbf{X}^t \boldsymbol{\epsilon}/\sqrt{n} \xrightarrow{\text{d}} N(\mathbf{0}, \mathbf{I})$. Therefore, because $0 \leq d_k^* \leq 1$, deduce that the second term on the right-hand side of (26) is $O_p(1)$. Now consider the first term on the right-hand side of (26). By our assumptions regarding the support of $\pi$ and $\mu$, we must have $d_k^* \geq \eta_0/(\eta_0 + s_0^2)$. Thus, $d_k^*$ remains bounded away from zero in probability. Hence, because $C_n/\sqrt{n} \to 0$, we have proven that

$$C_n^{-1} |\widehat{\beta}_{k,n}^*| \xrightarrow{\text{P}} \begin{cases} 0, & \text{if } k \in \mathscr{B}_0, \\ \infty, & \text{otherwise.} \end{cases} \quad \square$$

PROOF OF THEOREM 8. As the proof is somewhat lengthy, we first give a brief sketch. The basis for the proof will rely on the following result:

$$\widehat{\beta}_k^\circ[k] \xrightarrow{\text{P}} \begin{cases} 0, & \text{if } k_0 < k \leq K, \\ \beta_{k_0,0}, & \text{if } k = k_0, \\ \beta_{k,0} + \Delta_{k,0}, & \text{if } 1 \leq k < k_0, \end{cases} \tag{27}$$

where $\Delta_{k,0}$ is the $k$th coordinate of $\boldsymbol{\Sigma}_0[k:k]^{-1} \boldsymbol{\Sigma}_0[k:-k] \boldsymbol{\beta}_0[-k]$. Here $\boldsymbol{\beta}_0[-k] = (\beta_{0,k+1}, \ldots, \beta_{0,K})^t$, while $\boldsymbol{\Sigma}_0[k:k]$ and $\boldsymbol{\Sigma}_0[k:-k]$ are the $k \times k$ and $k \times (K-k)$ matrices associated with $\boldsymbol{\Sigma}_0$ which has been partitioned according to

$$\boldsymbol{\Sigma}_0 = \begin{pmatrix} \boldsymbol{\Sigma}_0[k:k] & \boldsymbol{\Sigma}_0[k:-k] \\ \boldsymbol{\Sigma}_0[-k:k] & \boldsymbol{\Sigma}_0[-k:-k] \end{pmatrix}.$$

First consider what (27) implies when $k < k_0$. Recall the definition (15) for $\widetilde{Z}_{k,n}$. Using $\widehat{\sigma}_n^2 \xrightarrow{\text{P}} \sigma_0^2$ (shown in the proof of Theorem 6) and that $s_{kk}[k]$ converges to the $k$th diagonal value of $\boldsymbol{\Sigma}_0[k:k]^{-1}$, a strictly positive value, deduce from the second limit of (27) that $\mathbb{P}\{|\widetilde{Z}_{k_0,n}| \geq z_{\alpha_{k_0}/2}\} \to 1$. Thus, for (a),

$$\mathbb{P}\{\hat{k}_B = k\} = \mathbb{P}\{|\widetilde{Z}_{k,n}| \geq z_{\alpha_k/2} \text{ and } |\widetilde{Z}_{j,n}| < z_{\alpha_j/2} \text{ for } j = k+1, \ldots, K\}$$
$$\leq \mathbb{P}\{|\widetilde{Z}_{k_0,n}| < z_{\alpha_{k_0}/2}\} \to 0.$$

For (b), observe that $\Delta_{k,0} = 0$ (by our assumption of orthogonality). Thus, when $k \leq k_0$, the last two lines of (27) imply that $\mathbb{P}\{|\widetilde{Z}_{k,n}| \geq z_{\alpha_k/2}\} \to 1$, and therefore,

$$\mathbb{P}\{\hat{k}_F = k - 1\} = \mathbb{P}\{|\widetilde{Z}_{k,n}| < z_{\alpha_k/2} \text{ and } |\widetilde{Z}_{j,n}| \geq z_{\alpha_j/2} \text{ for } j = 1, \ldots, k - 1\}$$
$$\leq \mathbb{P}\{|\widetilde{Z}_{k,n}| < z_{\alpha_k/2}\} \to 0.$$

Now for (c), due to orthogonality, $\widehat{Z}_{k,n} = \widetilde{Z}_{k,n}$ for $\widehat{Z}_{k,n}$ defined by (12). Thus,

$$\mathbb{P}\{\hat{k}_O \geq k_0\} \geq \mathbb{P}\{|\widetilde{Z}_{j,n}| \geq z_{\alpha_j/2} \text{ for } j = 1, \ldots, k_0\} \to 1.$$

Thus, for all three estimators the probability of the event $\{k < k_0\}$ tends to zero.

Now consider when $k > k_0$. We will show for (a) [and, therefore, for (b) and (c)]

(28) $$\widetilde{\mathbf{Z}}_n = (\widetilde{Z}_{k_0+1,n}, \ldots, \widetilde{Z}_{K,n})^t \stackrel{\mathrm{d}}{\rightsquigarrow} \mathrm{N}(\mathbf{0}_{K-k_0}, \mathbf{I}),$$

which implies $\{\widetilde{Z}_{k_0+1,n}, \ldots, \widetilde{Z}_{K,n}\}$ are asymptotically independent. By (28),

$$\mathbb{P}\{\hat{k}_B = k\} \to \mathbb{P}\{|\mathrm{N}(0,1)| \geq z_{\alpha_k/2}\} \prod_{j=k+1}^{K} \mathbb{P}\{|\mathrm{N}(0,1)| < z_{\alpha_j/2}\},$$

which is the third expression in (a). For (b), using (28) and the assumed orthogonality,

$$\mathbb{P}\{\hat{k}_F = k\} \to \mathbb{P}\{|\mathrm{N}(0,1)| < z_{\alpha_{k+1}/2}\} \prod_{j=k_0+1}^{k} \mathbb{P}\{|\mathrm{N}(0,1)| \geq z_{\alpha_j/2}\}.$$

Meanwhile, for OLS-hard, when $k > k_0$ or $k = k_0$,

$$\mathbb{P}\{\hat{k}_O = k\} \to \mathbb{P}\left\{\sum_{j=k_0+1}^{K} \mathbb{I}\{|Z_j| \geq z_{\alpha_j/2}\} = k - k_0\right\},$$

where $\{Z_{k_0+1}, \ldots, Z_K\}$ are mutually independent $\mathrm{N}(0,1)$ variables. This is the second expression in (c). Deduce that (a), (b) and (c) hold (the case $k = k_0$ for $k_B$ and $k_F$ can be worked out using similar arguments).

This completes the outline of the proof. Now we must prove (27) and (28). We start with (27). Let $\boldsymbol{\beta}_0[k] = (\beta_{0,1}, \ldots, \beta_{0,k})^t$. Some simple algebra shows that

(29) $$\widehat{\boldsymbol{\beta}}^\circ[k] = \boldsymbol{\beta}_0[k] + (\mathbf{X}[k]^t\mathbf{X}[k])^{-1}\mathbf{X}[k]^t\mathbf{X}[-k]\boldsymbol{\beta}_0[-k]$$
$$+ (\mathbf{X}[k]^t\mathbf{X}[k])^{-1}\mathbf{X}[k]^t\boldsymbol{\epsilon},$$

where $\mathbf{X}[-k]$ refers to the design matrix which excludes the first $k$ columns of $\mathbf{X}$. It is easy to show that the third term on the right-hand side is $o_p(1)$. Thus, it follows that

$$\widehat{\boldsymbol{\beta}}^\circ[k] \xrightarrow{\mathrm{p}} \boldsymbol{\beta}_0[k] + \boldsymbol{\Sigma}_0[k\!:\!k]^{-1}\boldsymbol{\Sigma}_0[k\!:\!-k]\boldsymbol{\beta}_0[-k],$$

which is what (27) asserts.

Finally, we prove (28). By (29), $\widehat{\beta}^\circ_k[k]$ is the $k$th coordinate of $(\mathbf{X}[k]^t\mathbf{X}[k])^{-1} \times \mathbf{X}[k]^t\boldsymbol{\epsilon}$ when $k > k_0$, and thus,

$$\widetilde{Z}_{k,n} = (\mathbf{0}^t_{k-1}, (s_{kk}[k])^{-1/2})(\mathbf{X}[k]^t\mathbf{X}[k]/n)^{-1}(\widehat{\sigma}_n^{-1}n^{-1/2}\mathbf{X}[k]^t\boldsymbol{\epsilon})$$
$$= (\widetilde{\mathbf{v}}^t_k, \mathbf{0}^t_{K-k})\boldsymbol{\xi}_n,$$

where $\boldsymbol{\xi}_n = \widehat{\sigma}_n^{-1}\mathbf{X}^t\boldsymbol{\epsilon}/\sqrt{n}$ and $\widetilde{\mathbf{v}}_k$ is the $k$-dimensional vector defined by

$$\widetilde{\mathbf{v}}^t_k = (\mathbf{0}^t_{k-1}, (s_{kk}[k])^{-1/2})(\mathbf{X}[k]^t\mathbf{X}[k]/n)^{-1}.$$

This allows us to write $\widetilde{\mathbf{Z}}_n = \mathbf{V}_n\boldsymbol{\xi}_n$, where

$$\mathbf{V}_n := \begin{bmatrix} \mathbf{v}^t_{k_0+1} \\ \vdots \\ \mathbf{v}^t_K \end{bmatrix} := \begin{bmatrix} (\widetilde{\mathbf{v}}^t_{k_0+1}, \mathbf{0}^t_{K-k_0-1}) \\ \vdots \\ \widetilde{\mathbf{v}}^t_K \end{bmatrix}.$$

Thus, because $\boldsymbol{\xi}_n \xrightarrow{\mathrm{d}} \mathrm{N}(\mathbf{0}, \boldsymbol{\Sigma}_0)$ by Lemma A.1, we have

$$\widetilde{\mathbf{Z}}_n \xrightarrow{\mathrm{d}} \mathrm{N}(\mathbf{0}_{K-k_0}, \mathbf{V}_0\boldsymbol{\Sigma}_0\mathbf{V}^t_0),$$

where $\mathbf{V}_0$ is the limit of $\mathbf{V}_n$. In particular, $\mathbf{V}_n\boldsymbol{\Sigma}_n\mathbf{V}^t_n \to \mathbf{V}_0\boldsymbol{\Sigma}_0\mathbf{V}^t_0$. To complete the proof, we show $\mathbf{V}_0\boldsymbol{\Sigma}_0\mathbf{V}^t_0 = \mathbf{I}$ by proving that $\mathbf{V}_n\boldsymbol{\Sigma}_n\mathbf{V}_n = \mathbf{I}$. Note by tedious (but straightforward) algebra that $\mathbf{v}^t_k\boldsymbol{\Sigma}_n\mathbf{v}_k = 1$. Consider $\mathbf{v}^t_j\boldsymbol{\Sigma}_n\mathbf{v}_k$ when $j \neq k$ and $j > k_0$. By (29), when $k > k_0$,

$$\widehat{\beta}^\circ_k[k] = n^{-1}s_{kk}[k]^{1/2}\mathbf{v}^t_k\mathbf{X}^t\boldsymbol{\epsilon}.$$

By Remark 8 we know that $\widehat{\beta}^\circ_{k_0+1}[k_0+1], \ldots, \widehat{\beta}^\circ_K[K]$ are uncorrelated. Thus $\mathbb{E}(\widehat{\beta}^\circ_j[j]\widehat{\beta}^\circ_k[k]) = 0$ if $j \neq k$, and therefore,

$$0 = \mathbb{E}(\mathbf{v}^t_j\mathbf{X}^t\boldsymbol{\epsilon}\mathbf{v}^t_k\mathbf{X}^t\boldsymbol{\epsilon}) = \mathbf{v}^t_j\mathbf{X}^t\mathbb{E}(\boldsymbol{\epsilon}\boldsymbol{\epsilon}^t)\mathbf{X}\mathbf{v}_k = \sigma^2_0\mathbf{v}^t_j\mathbf{X}^t\mathbf{X}\mathbf{v}_k.$$

Thus, $\mathbf{v}^t_j\boldsymbol{\Sigma}_n\mathbf{v}_k = 0$. Deduce that $\mathbf{V}_n\boldsymbol{\Sigma}_n\mathbf{V}^t_n = \mathbf{I}$ and, hence, that $\mathbf{V}_0\boldsymbol{\Sigma}_0\mathbf{V}^t_0 = \mathbf{I}$. □

**SVS Gibbs sampler.**

ALGORITHM. The SVS procedure uses a Gibbs sampler to simulate posterior values

$$(\boldsymbol{\beta}, \boldsymbol{J}, \boldsymbol{\tau}, w, \sigma^2 | \mathbf{Y}^*)$$

from (6), where $\boldsymbol{J} = (\mathscr{I}_1, \ldots, \mathscr{I}_K)^t$ and $\boldsymbol{\tau} = (\tau_1, \ldots, \tau_K)^t$. Recall that $\gamma_k = \mathscr{I}_k \tau_k^2$, so simulating $\boldsymbol{J}$ and $\boldsymbol{\tau}$ provides a value for $\boldsymbol{\gamma}$. The Gibbs sampler works as follows:

1. Simulate $(\boldsymbol{\beta} | \boldsymbol{\gamma}, \sigma^2, \mathbf{Y}^*) \sim \mathrm{N}(\boldsymbol{\mu}, \sigma^2 \boldsymbol{\Sigma})$, the conditional distribution for $\boldsymbol{\beta}$, where

$$\boldsymbol{\mu} = \boldsymbol{\Sigma} \mathbf{X}^t \mathbf{Y}^* \quad \text{and} \quad \boldsymbol{\Sigma} = (\mathbf{X}^t \mathbf{X} + \sigma^2 n \boldsymbol{\Gamma}^{-1})^{-1}.$$

2. Simulate $\mathscr{I}_k$ from its conditional distribution

$$(\mathscr{I}_k | \boldsymbol{\beta}, \boldsymbol{\tau}, w) \stackrel{\mathrm{ind}}{\sim} \frac{w_{1,k}}{w_{1,k} + w_{2,k}} \delta_{v_0}(\cdot) + \frac{w_{2,k}}{w_{1,k} + w_{2,k}} \delta_1(\cdot), \qquad k = 1, \ldots, K,$$

where

$$w_{1,k} = (1-w) v_0^{-1/2} \exp\left(-\frac{\beta_k^2}{2 v_0 \tau_k^2}\right)$$

and

$$w_{2,k} = w \exp\left(-\frac{\beta_k^2}{2 \tau_k^2}\right).$$

3. Simulate $\tau_k^{-2}$ from its conditional distribution,

$$(\tau_k^{-2} | \boldsymbol{\beta}, \boldsymbol{J}) \stackrel{\mathrm{ind}}{\sim} \mathrm{Gamma}\left(a_1 + \frac{1}{2}, a_2 + \frac{\beta_k^2}{2 \mathscr{I}_k}\right), \qquad k = 1, \ldots, K.$$

4. Simulate $w$, the complexity parameter, from its conditional distribution,

$$(w | \boldsymbol{\gamma}) \sim \mathrm{Beta}(1 + \#\{k : \gamma_k = 1\}, 1 + \#\{k : \gamma_k = v_0\}).$$

5. Simulate $\sigma^{-2}$ from its conditional distribution,

$$(\sigma^{-2} | \boldsymbol{\beta}, \mathbf{Y}^*) \sim \mathrm{Gamma}\left(b_1 + \frac{n}{2}, b_2 + \frac{1}{2n} \|\mathbf{Y}^* - \mathbf{X} \boldsymbol{\beta}\|^2\right).$$

6. This completes one iteration. Update $\boldsymbol{\gamma}$ by setting $\gamma_k = \mathscr{I}_k \tau_k^2$ for $k = 1, \ldots, K$.

COMPUTATIONS FOR LARGE $K$. The most costly computation in running the Gibbs sampler is the inversion

$$\boldsymbol{\Sigma} = (\mathbf{X}^t \mathbf{X} + \sigma^2 n \boldsymbol{\Gamma}^{-1})^{-1}$$

required in updating $\boldsymbol{\beta}$ in step 1. This requires $\mathrm{O}(K^3)$ operations and can be tremendously slow when $K$ is large.

A better approach is to update $\boldsymbol{\beta}$ in $B$ blocks of size $q$. This will reduce computations to order $\mathrm{O}(B^{-2}K^3)$, where $K = Bq$. To proceed, decompose $\boldsymbol{\beta}$ as $(\boldsymbol{\beta}_{(1)}^t, \ldots, \boldsymbol{\beta}_{(B)}^t)^t$, $\boldsymbol{\Gamma}$ as $\mathrm{diag}(\boldsymbol{\Gamma}_{(1)}, \ldots, \boldsymbol{\Gamma}_{(B)})$ and $\mathbf{X}$ as $[\mathbf{X}_{(1)}, \ldots, \mathbf{X}_{(B)}]$. Now update each component $\boldsymbol{\beta}_{(j)}$, $j = 1, \ldots, B$, conditioned on the remaining values. Using a subscript $-(j)$ to indicate exclusion of the $j$th component, draw $\boldsymbol{\beta}_{(j)}$ from a $\mathrm{N}(\boldsymbol{\mu}_j, \sigma^2 \boldsymbol{\Sigma}_j)$ distribution, where

$$\boldsymbol{\mu}_j = \boldsymbol{\Sigma}_j \mathbf{X}_{(j)}^t (\mathbf{Y}^* - \mathbf{X}_{-(j)} \boldsymbol{\beta}_{-(j)}) \quad \text{and} \quad \boldsymbol{\Sigma}_j = (\mathbf{X}_{(j)}^t \mathbf{X}_{(j)} + \sigma^2 n \boldsymbol{\Gamma}_{(j)}^{-1})^{-1}.$$

Notice that the cross-product terms $\mathbf{X}_{(j)}^t \mathbf{X}_{(j)}$ and $\mathbf{X}_{(j)}^t \mathbf{X}_{-(j)}$ can be extracted from $\mathbf{X}^t \mathbf{X}$ and do not need to be computed.

**Acknowledgments.** The authors thank Benedikt Pötscher and Hannes Leeb for helpful discussion surrounding Theorem 8.

## REFERENCES

Barbieri, M. and Berger, J. (2004). Optimal predictive model selection. *Ann. Statist* **32** 870–897. MR2065192

Bickel, P. and Zhang, P. (1992). Variable selection in non-parametric regression with categorical covariates. *J. Amer. Statist. Assoc.* **87** 90–97. MR1158628

Breiman, L. (1992). The little bootstrap and other methods for dimensionality selection in regression: $X$-fixed prediction error. *J. Amer. Statist. Assoc.* **87** 738–754. MR1185196

Chipman, H. (1996). Bayesian variable selection with related predictors. *Canad. J. Statist.* **24** 17–36. MR1394738

Chipman, H. A., George, E. I. and McCulloch, R. E. (2001). The practical implementation of Bayesian model selection (with discussion). In *Model Selection* (P. Lahiri, ed.) 65–134. IMS, Beachwood, OH. MR2000752

Clyde, M., DeSimone, H. and Parmigiani, G. (1996). Prediction via orthogonalized model mixing. *J. Amer. Statist. Assoc.* **91** 1197–1208.

Clyde, M., Parmigiani, G. and Vidakovic, B. (1998). Multiple shrinkage and subset selection in wavelets. *Biometrika* **85** 391–401. MR1649120

Efron, B., Hastie, T., Johnstone, I. and Tibshirani, R. (2004). Least angle regression (with discussion). *Ann. Statist.* **32** 407–499. MR2060166

George, E. I. (1986). Minimax multiple shrinkage estimation. *Ann. Statist.* **14** 188–205. MR829562

George, E. I. and McCulloch, R. E. (1993). Variable selection via Gibbs sampling. *J. Amer. Statist. Assoc.* **88** 881–889.

Geweke, J. (1996). Variable selection and model comparison in regression. In *Bayesian Statistics 5* (J. M. Bernardo, J. O. Berger, A. P. Dawid and A. F. M. Smith, eds.) 609–620. Oxford Univ. Press, New York. MR1425430

Hoerl, A. E. (1962). Application of ridge analysis to regression problems *Chemical Engineering Progress* **58** 54–59.

Hoerl, A. E. and Kennard, R. W. (1970). Ridge regression: Biased estimation for nonorthogonal problems. *Technometrics* **12** 55–67.

Ishwaran, H. (2004). Discussion of "Least angle regression," by B. Efron, T. Hastie, I. Johnstone and R. Tibshirani. *Ann. Statist.* **32** 452–457.

Ishwaran, H. and Rao, J. S. (2000). Bayesian nonparametric MCMC for large variable selection problems. Unpublished manuscript.

Ishwaran, H. and Rao, J. S. (2003). Detecting differentially expressed genes in microarrays using Bayesian model selection *J. Amer. Statist. Assoc.* **98** 438–455. MR1995720

Ishwaran, H. and Rao, J. S. (2005). Spike and slab gene selection for multigroup microarray data. *J. Amer. Statist. Assoc.* To appear. MR1995720

Knight, K. and Fu, W. (2000). Asymptotics for lasso-type estimators *Ann. Statist.* **28** 1356–1378. MR1805787

Kuo, L. and Mallick, B. K. (1998). Variable selection for regression models. *Sankhyā Ser. B* **60** 65–81. MR1717076

Le Cam, L. and Yang, G. L. (1990). *Asymptotics in Statistics*: *Some Basic Concepts*. Springer, New York. MR1066869

Leeb, H. and Pötscher, B. M. (2003). The finite-sample distribution of post-model-selection estimators, and uniform versus non-uniform approximations. *Econometric Theory* **19** 100–142. MR1965844

Lempers, F. B. (1971). *Posterior Probabilities of Alternative Linear Models*. Rotterdam Univ. Press. MR359234

Mitchell, T. J. and Beauchamp, J. J. (1988). Bayesian variable selection in linear regression (with discussion). *J. Amer. Statist. Assoc.* **83** 1023–1036. MR997578

Pötscher, B. M. (1991). Effects of model selection on inference. *Econometric Theory* **7** 163–185. MR1128410

Rao, C. R. and Wu, Y. (1989). A strongly consistent procedure for model selection in a regression problem. *Biometrika* **76** 369–374. MR1016028

Rao, J. S. (1999). Bootstrap choice of cost complexity for better subset selection. *Statist. Sinica* **9** 273–287. MR1678894

Shao, J. (1993). Linear model selection by cross-validation *J. Amer. Statist. Assoc.* **88** 486–494. MR1224373

Shao, J. (1996). Bootstrap model selection. *J. Amer. Statist. Assoc.* **91** 655–665. MR1395733

Shao, J. (1997). An asymptotic theory for linear model selection (with discussion). *Statist. Sinica* **7** 221–264. MR1466682

Shao, J. and Rao, J. S. (2000). The GIC for model selection: A hypothesis testing approach. Linear models. *J. Statist. Plann. Inference* **88** 215–231. MR1792042

Zhang, P. (1992). On the distributional properties of model selection criteria. *J. Amer. Statist. Assoc.* **87** 732–737. MR1185195

Zhang, P. (1993). Model selection via multifold cross validation. *Ann. Statist.* **21** 299–313. MR1212178

Zheng, X. and Loh, W.-Y. (1995). Consistent variable selection in linear models. *J. Amer. Statist. Assoc.* **90** 151–156. MR1325122

Zheng, X. and Loh, W.-Y. (1997). A consistent variable selection criterion for linear models with high-dimensional covariates. *Statist. Sinica* **7** 311–325. MR1466685

Department of Biostatistics
and Epidemiology
Wb4
Cleveland Clinic Foundation
9500 Euclid Avenue
Cleveland, Ohio 44195
USA
e-mail: ishwaran@bio.ri.ccf.org

Department of Epidemiology
and Biostatistics
Case Western Reserve University
10900 Euclid Avenue
Cleveland, Ohio 44106
USA
e-mail: sunil@hal.epbi.cwru.edu